\documentclass[a4paper,11pt]{article}

\usepackage{amssymb}
\usepackage{latexsym}
\usepackage{amsmath,amsthm}
\usepackage{xy}
\usepackage{xypic}

\numberwithin{equation}{section}

\newtheorem{thm}{Theorem}[section]
\newtheorem{cor}[thm]{Corollary}
\newtheorem{lem}[thm]{Lemma}
\newtheorem{prop}[thm]{Proposition}
\newtheorem{defn}[thm]{Definition}
\newtheorem{rem}[thm]{Remark}

\def\CN{{\mathcal {N}}}

\def\CB{{\mathcal {B}}}

\def\CH{{\mathcal {H}}}
\def\CK{{\mathcal {K}}}

\def\CD{{\mathcal {D}}}
\def\CO{{\mathcal {O}}}

\def\CR{{\mathcal {R}}}

\def\ua{{\underline{a}}}

\def\uT{{\underline{T}}}
\def\uR{{\underline{R}}}
\def\uV{{\underline{V}}}

\def\uA{{\underline{A}}}
\def\uB{{\underline{B}}}
\def\uC{{\underline{C}}}

\def\uE{{\underline{E}}}

\def\uY{{\underline{Y}}}

\def\uR{{\underline{R}}}

\def\ua{{\underline{a}}}
\def\ub{{\underline{b}}}
\def\uc{{\underline{c}}}

\def\N{{\mathbb {N}}}
\def\C{{\mathbb {C}}}

\newcommand{\eins}{\mathbf{1}}

\newcommand{\Od}{\mathcal{O}_d}

\newcommand{\DsA}{D_{*,A}}

\begin{document}

\begin{center}
{\Large {\bf Characteristic Functions of Liftings}}
\end{center}

\vspace{1cm}

\begin{center}Santanu Dey and Rolf Gohm
\end{center}

\vspace{1cm}

{\bf Abstract:} {\em We introduce characteristic functions for certain
contractive liftings of row contractions. These are multi-analytic
operators which classify the liftings up to unitary equivalence and
provide a kind of functional model. The most important cases are
subisometric and coisometric liftings. We also identify the most
general setting which we call reduced liftings. We derive properties
of these new characteristic functions and discuss the relation to
Popescu's definition for completely non-coisometric row
contractions. Finally we apply our theory to completely positive
maps and prove a one-to-one correspondence between the fixed point
sets of completely positive maps related to each other by a
subisometric lifting.}

\vspace{2cm}
\noindent {\bf keywords:} characteristic function, contractive lifting,  row
contraction, multi-analytic, subisometric dilation,
coisometric, completely non-coisometric, completely positive
map, fixed point
\vspace{2cm}

\section*{Introduction}

Let $C$ be a contraction on a Hilbert space $\CH_C$. Then a contraction $E$ on a Hilbert space
$\CH_E \supset \CH_C$ is called a contractive lifting of $C$ if $P E = C P$, where $P$ is the orthogonal projection from $\CH_E$ onto $\CH_C$. In other words, we have an operator matrix
\begin{align} \label{0.1}
E = \left(
\begin{array}{cc}
C & 0 \\
B & A \\
\end{array}
\right).
\end{align}
See Chapter 5 of \cite{FF90}. In this book C.\,Foias and A.E.\,Frazho amply demonstrate the importance of understanding the structure of contractive liftings, in particular in connection with the commutant lifting theorem and its applications.

The minimal isometric dilation (mid for short) of $C$ is the most prominent example of a contractive lifting. In \cite{DF84} R.G.\,Douglas and C.\,Foias introduced subisometric dilations (see also
Chapter 8.3 of \cite{Ber88} for a discussion closer to our point of view).
These are contractive liftings with the property that the mid of $E$ is also minimal as an isometric dilation of $C$. In this context Douglas and Foias were especially interested in problems of uniqueness and of commutant lifting. We arrived at the subisometric property in a completely different way and ask different questions about it. Let us briefly describe the most relevant aspects of this development.

Many results of the Sz.-Nagy/Foias-theory for contractions \cite{NF70} can be generalized to row contractions $\underline{C} = (C_1,\ldots,C_d)$, i.e. tuples of operators such that
$\sum^d_{i=1} C_i C^*_i \leq \eins$. This has been done very systematically by G.\,Popescu starting
with \cite{Po89a} and many people contributed to this development, an incomplete list of work related to our interests is \cite{Ar98,BBD04,BDZ06,BES05,DKS01,Po89b,Po03,Po05}.
In particular in \cite{Po89b} G.\,Popescu described a class of
multi-analytic operators which classify completely non-coisometric (c.n.c.) row contractions up to
unitary equivalence and called them characteristic functions, in analogy to a similar concept in the
Sz.-Nagy/Foias-theory. In \cite{DG07} S.\,Dey and R.\,Gohm started from some seemingly unrelated questions in noncommutative probability theory arising in \cite{Go04,GKL06} and established a class of multi-analytic operators which are associated to certain rather special coisometric row contractions
(i.e., $\sum^d_{i=1} C_i C^*_i = \eins$). Investigating their properties we came to the conclusion that there are good reasons to think of them as of characteristic functions for these tuples. This is not covered by Popescu's theory.

In this paper we will show that it is the property of being a subisometric lifting which makes this analysis possible. This is a vast generalization of the setting of \cite{DG07} and it clarifies the mechanism behind it. It is straightforward to define liftings for row contractions. Let
$\uE = (E_1,\ldots,E_d)$ be a row contraction on a Hilbert space $\CH_E \supset \CH_C$. If for all $i=1,\ldots,d$
(with $d$ countable) we have an operator matrix
\begin{align} \label{0.2}
E_i =\left(
\begin{array}{cc}
C_i  & 0 \\
B_i  & A_i \\
\end{array}
\right)
\end{align}
with respect to $\CH_C \oplus \CH_C^\perp$
then we say that $\uE$ is a lifting of $\uC=(C_1,\ldots,C_d)$ by $\uA=(A_1,\ldots,A_d)$ (or that $\uE$ is an extension of $\uA$
by $\uC$). The subisometric property in the form given here also makes sense for row contractions,
using Popescu's theory of mid for row contractions \cite{Po89a}. This is worked out in Section 1 below. It then turns out that there is a Beurling-type classification of subisometric liftings,
involving a correspondence to certain multi-analytic inner operators (Theorem 1.6). They classify subisometric liftings up to unitary equivalence, so we call them characteristic functions of
(subisometric) liftings.

In Section 2 we focus on coisometric liftings, i.e. $\sum^d_{i=1} E_i E^*_i = \eins$, emphasizing another type of classification which uses an isometry $\gamma$ mapping the defect space $\CD_{*,A}$ of $\uA$
into the defect space $\CD_C$ of $\uC$ (Theorem 2.1). The connection to Section 1 lies in the fact that coisometric liftings by $*-$stable $\uA$ are subisometric (Proposition 2.3). But this is only a special case and we have to generalize further.

This is done in Section 3. We get a hint from a result about contractive liftings for single contractions. Lemma 2.1 in Chap.IV of \cite{FF90} states that
$E = \left(
\begin{array}{cc}
C & 0 \\
B & A \\
\end{array}
\right)$
is a contraction if and only if $C$ and $A$ are contractions and there exists a contraction
$\gamma: \CD_{*,A} \rightarrow \CD_C$ such that
\begin{align} \label{0.3}
B = D_{*,A}\, \gamma^* \, D_C,
\end{align}
where $D_{*,A}$
and $D_C$ are the defect operators of $A^*$ and $C$. We establish an analogous result for
row contractions (Proposition 3.1). This shows that the isometry $\gamma$ occurring for
coisometric liftings in Section 2 has to be replaced in a more general setting by a contraction.

The most general situation where we can establish a satisfactory theory of characteristic functions for liftings is identified in Section 3 and we call such liftings reduced. The technical tool here
is to use the Wold decomposition for the mid's. For $\gamma$ we isolate the special property needed and call it resolving. Reduced liftings include subisometric liftings as well as coisometric liftings by c.n.c. row contractions. We define characteristic functions for reduced liftings (Definition 3.6)
and we argue that this is the most general setting which is natural for that. These characteristic functions are multi-analytic operators (not inner in general) and they characterize reduced liftings
up to unitary equivalence. They also provide a kind of functional model for the lifting which is useful for a closer investigation of the structure of the lifting in the same sense as the characteristic functions of Sz.-Nagy/Foias and of Popescu are useful in their context.

In Section 4 we study some further properties of these characteristic functions. In particular we clarify the connection to Popescu's characteristic functions and we investigate iterated liftings,
showing a factorization result for our characteristic functions (Theorem 4.1). This is another indication that our definition leads to a promising theory.

We believe that in particular the theory of subisometric liftings may be even more interesting for row contractions than it is for single contractions. There is a straightforward way to transfer results from a row contraction $\uC = (C_1,\ldots,C_d)$ to the completely positive map
$\Phi_C: X \mapsto \sum^d_{i=1} C_i X C^*_i$. This topic is taken up in Section 5. We define
characteristic functions for liftings of completely positive maps and show in which way they are characteristic in this case (Corollary 5.2). We investigate what subisometric lifting means in this context and prove a one-to-one correspondence between the fixed point sets
(Theorem 5.4). In particular
we consider the situation where a normal invariant state is restricted to its support (Corollary 5.6). From our point of view these applications give a strong motivation for further developing the theory of liftings for row contractions.

In an Appendix we reprove a commutant lifting theorem by O.\,Bratteli, P.\,Jorgensen, A.\,Kishimoto and R.F.\,Werner \cite{BJKW00}, used in Section 5, in a way that helps to understand its role in our theory.

\section{Subisometric Liftings}

In this section we define subisometric liftings in the setting of row contractions and show that there is a nice Beurling-type classification for them.
\\

We recall the notion of a minimal isometric dilation for a row contraction, cf. \cite{Po89a}.
Let $\uT=(T_1,\cdots,T_d)$ be a row contraction on a Hilbert space $\CH$.
Treating $\uT$ as an operator from $ \bigoplus^d_{i=1} \CH$
to $\CH$, define
$D_*:= (\eins-\uT \uT^*)^\frac{1}{2}: \CH \rightarrow \CH$ and
$D:=(\eins-\uT^*\uT)^\frac{1}{2}:
\bigoplus^d_{i=1} \CH \rightarrow \bigoplus^d_{i=1} \CH$.
This implies that
\begin{align}\label{1.1}
D_*= (\eins -\sum^d_{i=1} T_i T^*_i)^\frac{1}{2},
~~~D=(\delta_{ij} \eins -T_i^*T_j)^\frac{1}{2}_{d\times d}
\end{align}
Let $\CD_*:= \overline{\mbox{range~} D_*}$ and $\CD:=\overline{\mbox{range~} D}$.

We  use the following {\em multi-index notation}.
Let $\Lambda $ denote the set $\{ 1, 2, \ldots , d\}$ and
$\tilde{\Lambda}:=\cup_{n=0}^{\infty} \Lambda^n$, where $\Lambda ^0:=\{ 0\}$.
If $\alpha \in \Lambda^n \subset \tilde{\Lambda}$ the integer $n = |\alpha|$
is called its length.
Now $T_\alpha$ with $\alpha = (\alpha_1, \cdots ,\alpha_n) \in \Lambda^n$ means
$T_{\alpha_1} T_{\alpha_2} \ldots T_{\alpha_n}$.

The full Fock space over $\C^d$ ($d\geq 2$) denoted by $\Gamma (\C^d)$ is
\begin{align}\label{1.2}
\Gamma (\C^d) := \mathbb{C}\oplus \C^d \oplus (\C^d)^{\otimes ^2}\oplus \cdots
\oplus (\C^d) ^{\otimes ^m}\oplus \cdots.
\end{align}
To simplify notation we shall often only write $\Gamma$ instead of $\Gamma (\C^d)$. The vector
$e_0 := 1\oplus 0\oplus \cdots$ is called the vacuum vector. Let
$e_1, \ldots , e_d$ be the standard
orthonormal basis of $\C^d$. We include $d=\infty$
in which case $\C^d$ stands for a complex separable Hilbert space
of infinite dimension.
For $\alpha \in \Lambda^n$, $e_{\alpha}$ will denote the vector $e_{\alpha
_1}\otimes e_{\alpha _2}\otimes \cdots \otimes e_{\alpha _n}$ in
the full Fock space $\Gamma$. Then $e_{\alpha}$ over all $\alpha \in \tilde{\Lambda}$
forms an orthonormal basis of the full Fock space.
The (left) creation operators
$L_i$ on $\Gamma({\mathbb{C}^d})$ are defined by
$ L_i x = e_i \otimes x $  for $1 \leq i \leq d$ and $x
\in \Gamma({\mathbb{C}}^d).$ Then
$\underline{L}= (L_1, \ldots , L_d)$ is a row isometry, i.e., the $L_i$ are isometries with
orthogonal ranges.
\\

Using the definition of lifting in the introduction
a minimal isometric dilation (mid for short) can be described as an isometric lifting
$\uV$ of $\uT$ such that the spaces
$V_\alpha \CH$ with $\alpha \in \tilde{\Lambda}$ together span the Hilbert space
on which the $V_i$ are defined. It is an important fact, which we shall
use repeatedly, that such minimal isometric dilations are unique up to unitary
equivalence (cf. \cite{Po89a}).
A useful model for the mid is given by a version of the Sch\"{a}ffer construction,
given in \cite{Po89a}. Namely, we can realize a mid $\uV$ of $\uT$
on the Hilbert space $\hat{\CH} := \CH  \oplus (\Gamma \otimes \CD)$,
\begin{align}\label{1.3}
V_i(h \oplus \sum_{\alpha \in \tilde{\Lambda}}  e_\alpha \otimes d_\alpha)
= T_i h \oplus [e_0 \otimes D_i h + e_i \otimes \sum_{\alpha \in
\tilde{\Lambda}} e_\alpha \otimes d_\alpha]
\end{align}
for $h \in \CH$ and $d_\alpha \in \CD.$ Here
$D_i h := D (0,\ldots,0,h,0,\ldots,0)$ and $h$ is embedded at the
$i^{th}$ component.

If we have more than one row contraction at the same time
then we shall use the above notations with superscripts or subscripts, as convenient.
We are now ready for the basic definition in this section.
\begin{defn}
Let $\uC=(C_1,\cdots,C_d)$ be a row contraction on a Hilbert space $\CH_C$.
A lifting $\underline{E}$ of $\underline{C}$ on $\CH_E \supset \CH_C$ is called
subisometric if the corresponding mids $\underline{V}^E$
(on the Hilbert space $\hat{\CH}_E$) and $\underline{V}^C$ (on the Hilbert space $\hat{\CH}_C$) are unitarily equivalent,
in the sense that there exists a unitary $W: \hat{\CH}_E \rightarrow \hat{\CH}_C$ such that
$W |_{\CH_C} = \eins |_{\CH_C}$ and $W V_i^E = V_i^C W$.
\end{defn}
For $d=1$ this is consistent with the definition of subisometric dilation in \cite{DF84}, see the discussion in the introduction.
Note that the mid $\underline{V}^C$ is an example of a subisometric lifting in this sense. Another
(trivial) example is $\uC$ itself (considered as a lifting of $\uC$).
Further note that, given the mids $\underline{V}^E$ and $\underline{V}^C$,
the unitary $W$ is uniquely determined by its properties (use the minimality of $\underline{V}^C$).
\\

We want to make the structure of subisometric liftings more explicit. Let
$\underline{E} = (E_1,\ldots,E_d)$ be a subisometric lifting of $\uC=(C_1,\cdots,C_d)$
on ${\CH_E} = {\CH}_C \oplus {\CH}_A$ as in Definition 1.1, so that
for all $i=1,\ldots,d$ we have block matrices
\begin{align}\label{1.4}
E_i =\left(
\begin{array}{cc}
C_i  & 0 \\
B_i  & A_i \\
\end{array}
\right)
\end{align}
Let $\underline{V}^C$ be the mid of $\underline{C}$, realized as in \eqref{1.3} on the space
$\hat{\CH}_C = \CH_C  \oplus (\Gamma \otimes \CD_C)$. Because ${\CH_E} = {\CH}_C \oplus {\CH}_A \subset \hat{\CH}_E$ we can use the unitary $W$ from the subisometric lifting property to obtain a subspace $\CH_{A*} := W \CH_A \subset \Gamma \otimes \CD_C$. Further
$\CH_{E*} := {\CH}_C \oplus {\CH}_{A*} \subset \hat{\CH}_C$, and $\underline{V}^C$ is also a mid
of the row contraction $\uE_*$ which is transferred by $W$ from the unitarily equivalent original $\uE$.
We can write
\begin{align}\label{1.5}
E_{i*} =\left(
\begin{array}{cc}
C_i  & 0 \\
B_{i*} & A_{i*} \\
\end{array}
\right)
\end{align}
so $\uE_*$ is also a lifting of $\uC$.

Because $\underline{V}^C$ is a mid of $\uE_*$ it follows that $\CH_{E*}$ is coinvariant for
$\uV^C$ (by which we mean that it is invariant for all $(V^C_i)^*$, $i=1,\ldots,d$).
Note that
\begin{align}\label{1.6}
V_i^C \;|_{\Gamma \otimes \CD_C} = L_i \otimes \eins.
\end{align}
Hence $\underline{L \otimes \eins}$ is an isometric lifting of $\uA_*$,
in particular $\CH_{A*}$ is
coinvariant for $\underline{L \otimes \eins}$. An isometric lifting always contains the mid.
In particular the
mid of $\uA_*$ lives on the space $\overline{span}\{(L_\alpha \otimes \eins) \CH_{A*},
\alpha \in \tilde{\Lambda} \}$. This subspace is reducing for the
$L_i \otimes \eins$ for all $i=1,\ldots,d$ and hence has the form
$\Gamma \otimes {\cal E}$ for a subspace ${\cal E}$ of $\CD_C$,
see for example Cor.1.7 of \cite{Po05}, where it is done in a
more general setting. In this reference the space ${\cal E}$
is described as the closure of the image of $\CH_{A*}$ under the orthogonal
projection onto $e_0 \otimes \CD_C$.

We can obtain a more concrete formula for ${\cal E}$ by comparing this result with another way of writing the mid.
First note that, as a compression of $\underline{L \otimes \eins}$, the row contraction
$\uA_*$ (and hence also $\uA$) is $*$-stable, i.e., for all $h \in \CH_A$
\begin{align}\label{1.7}
\lim_{n\to\infty} \sum_{|\alpha|=n} \| A^*_\alpha h\|^2 = 0\;,
\end{align}
cf. \cite{Po89a}, Prop.2.3 (where it is called pure).
In this case, with $D_{*,A} \;=\; (\eins- \uA \uA^*)^{\frac{1}{2}}:\; \CH_A \rightarrow \CH_A$ and $\CD_{*,A}$ its closed range, the map
\begin{align}\label{1.8}
\CH_A \to \Gamma \otimes \CD_{*,A}
\\
h \mapsto \sum_{\alpha \in \tilde{\Lambda}} e_\alpha \otimes D_{*,A} A^*_\alpha h \nonumber
\end{align}
is isometric (Popescu's Poisson kernel, cf. \cite{Po03}). With this embedding of $\CH_A$
it can be checked that now $\underline{L \otimes \eins}$ on $\Gamma \otimes \CD_{*,A}$ is a mid of $\uA$.

Because mids are unique up to unitary equivalence we have a unitary $u: \Gamma \otimes \CD_{*,A} \rightarrow \Gamma \otimes {\cal E}$ such that $u \CH_A = \CH_{A*}$ and
$u (L_i \otimes \eins) = (L_i \otimes \eins) u$ for all $i=1,\ldots,d$. The commutation relation implies that $u$ is of the form $\eins \otimes u'$, where $u'$ is a unitary from
$\CD_{*,A}$ onto ${\cal E}$ (you may use the fact that $e_0 \otimes \CD_{*,A}$ respectively
$e_0 \otimes {\cal E}$ are the uniquely determined wandering subspaces).
Thinking of $u'$ as an isometry from $\CD_{*,A}$ into $\CD_C$ we call it $\gamma$.
So $\gamma: \CD_{*,A} \rightarrow \CD_C$ has ${\cal E}$ as its range and it is canonically
associated to a subisometric lifting in the way shown above.

Using $\gamma$ we see that the embedding of $\CH_{A}$ into $\Gamma \otimes {\CD_C}$ is  automatically of Poisson kernel type \eqref{1.8}, namely
\begin{align}\label{1.9}
\CH_A \ni h \mapsto \sum_{\alpha \in \tilde{\Lambda}}
e_\alpha \otimes \gamma D_{*,A} A^*_\alpha h \in \Gamma \otimes {\CD_C}
\end{align}
which is an explicit formula for
the embedding $W |_{\CH_A}: \CH_A \rightarrow \CH_{A*} \subset \Gamma \otimes \CD_C$.

Note also that the isometry $\gamma$ is closely related to the $\uB$-part of the lifting $\uE$.
In fact, because $E^*_{i*} = (V^C_i)^* |_{\CH_{E*}}$ we obtain $B^*_{i*} = p_C (V^C_i)^* p_{A*}$,
where $p_C, p_{A*}$ are the orthogonal projections onto $\CH_C, \CH_{A*}$. Combining
this with \eqref{1.3} and \eqref{1.9}
yields $B^*_i = D^*_{i,C}\, \gamma\, D_{*,A}: \CH_A \rightarrow \CH_C,\; i=1,\ldots,d$.
Or in a more compact form
\begin{align}\label{1.10}
\uB^* = D^*_C \, \gamma \, D_{*,A}.
\end{align}

\begin{prop}
A lifting $\uE$ of a row contraction $\uC$ with
\[
E_i =\left(
\begin{array}{cc}
C_i  & 0 \\
B_i  & A_i \\
\end{array}
\right),
\qquad i=1,\ldots,d,
\]
is subisometric if and only if $\uA$ is $*$-stable and $\uB = D_{*,A} \gamma^* D_C$
with an isometry $\gamma: \CD_{*,A} \rightarrow \CD_C$.
\end{prop}

\begin{proof}
We have already seen above that if $\uE$ is subisometric then the conditions are satisfied.
Conversely, if $\uA$ is $*$-stable then use the isometry $\gamma$ to embed $\uA$
(as $\uA_*$) and its mid into $\Gamma \otimes \CD_C$ as in \eqref{1.9}. Then the
formula for $\uB$ (or \eqref{1.10}) combined with \eqref{1.3} for $\uC$ shows that $\uV^C$ is a mid for $\uE_*$ which is unitarily equivalent to $\uE$. (Clearly $\uV^C$ is minimal for
$\uE_*$ because it is already minimal for $\uC$.) Hence $\uE$ is subisometric.
\end{proof}

\begin{rem}
This is consistent with the results for $d=1$ in \cite{DF84} which we mentioned in the introduction.
$\gamma$ unitary corresponds to what Douglas and Foias call a minimal subisometric dilation. We have no reason for imposing this condition and continue to consider general subisometric liftings.
Compare also Chapter 8.3 of \cite{Ber88}.
\end{rem}

Classifying subisometric liftings becomes especially transparent by focusing on the invariant subspace associated to it.

\begin{defn}
Let $\uE$ on ${\CH_E} = {\CH}_C \oplus {\CH}_A$ be a subisometric lifting of $\uC$ on
${\CH}_C$, notation as in Definition 1.1. Then we call
\begin{align}\label{1.11}
{\cal N} := (\Gamma \otimes \CD_C) \ominus W \CH_A
\end{align}
the invariant subspace associated to the subisometric lifting.
Clearly ${\cal N}$ is invariant for $L_i \otimes \eins$, $i=1,\ldots,d$.
\end{defn}

We can go the way back. Let $\uC$ on ${\cal H}_C$ be a row contraction.
If ${\cal N} \subset \Gamma \otimes \CD_C$ is a subspace which is invariant for all
$L_i \otimes \eins, \; i=1,\ldots,d$ then we can define
\begin{align}\label{1.12,1.13}
\CH_{A*} := (\Gamma \otimes \CD_C) \ominus \CN
\\
\CH_* := \CH_C \oplus \CH_{A*}
\end{align}
On $\CH_C \oplus (\Gamma \otimes \CD_C)$ we have the mid $\uV^C$ of $\uC$, as in \eqref{1.3},
so we can further define
\begin{align}\label{1.14}
\uE_* = (E_{1*},\ldots,E_{d*}),\quad E_{i*} := P_{\CH_*} V_i^C |_{\CH_{E*}}:
\CH_{E*} \rightarrow \CH_{E*}
\end{align}
Then $\uE_*$ is a row contraction and
\begin{align}\label{1.15}
E_{i*} =\left(
\begin{array}{cc}
C_i  & 0 \\
B_{i*}  & A_{i*} \\
\end{array}
\right)
\end{align}
with respect to the decomposition $\CH_{E*} := \CH_C \oplus \CH_{A*}$, i.e., $\uE_*$ is a
lifting of $\uC$. Then $\uV^C$ is a mid of $\uE_*$ (minimal
because it is already minimal for $\uC$). Hence we have constructed a subisometric lifting.
We are back in the setting of Proposition 1.2.
\\

These considerations suggest a classification of subisometric liftings along a Beurling type theorem for the associated
invariant subspaces. It is instructive to introduce the generalized inner functions occurring here directly from the definition of subisometric lifting.

So let $\uE$ be a subisometric lifting of $\uC$.
Then the mids $\uV^E$ of $\uE$ and $\uV^C$ of $\uC$ are connected by the unitary
\begin{align}\label{1.16}
W: \hat{\CH}_E = \CH_E \oplus (\Gamma \otimes \CD_E) \rightarrow
\hat{\CH}_C = \CH_C \oplus (\Gamma \otimes \CD_C)
\end{align}
such that $W |_{\CH_C} = \eins |_{\CH_C}$ and $W V^E_i = V_i^C W$ for $i=1,\ldots,d$.
If we define the isometry
\begin{align}\label{1.17}
M_{C,E} := W |_{\Gamma \otimes \CD_E}
\end{align}
then from \eqref{1.3} and \eqref{1.16} we obtain
\begin{align}\label{1.18}
M_{C,E} (L_i \otimes \eins_E) = (L_i \otimes \eins_C) M_{C,E}
\end{align}
which means that
$M_{C,E}: \Gamma \otimes \CD_E \rightarrow \Gamma \otimes \CD_C$
is a {\em multi-analytic inner operator} determined by its {\em symbol}
\begin{align}\label{1.19}
\Theta_{C,E}: \CD_E \rightarrow \Gamma \otimes \CD_C,
\quad \Theta_{C,E} = W |_{e_0 \otimes \CD_E}.
\end{align}
according to the terminology introduced in \cite{Po89b}.
Obviously this is nothing but the multi-analytic inner operator corresponding to the invariant subspace $\CN$, in fact it is easy to check that
\begin{align}\label{1.20}
\CN = M_{C,E} (\Gamma \otimes \CD_E),
\end{align}
compare the Beurling type theorem in \cite{Po89b}.
Our new insight is that it is connected to the subisometric lifting $\uE$ of $\uC$.

\begin{defn}
We call $M_{C,E}$ (or $\Theta_{C,E}$) the characteristic function of the
subisometric lifting $\uE$ of $\uC$.
\end{defn}

It is not difficult to check that two multi-analytic inner operators
$M: \Gamma \otimes \CD \rightarrow \Gamma \otimes {\cal E}$ and
$M^\prime: \Gamma \otimes \CD^\prime \rightarrow \Gamma \otimes {\cal E}$
with symbols $\Theta, \Theta^\prime$ describe the same invariant subspace
if and only if there exists a unitary $v: \CD \rightarrow \CD^\prime$ such that
$\Theta = \Theta^\prime v$. Let us call multi-analytic functions {\em equivalent}
if they are related in this way. We are ready for our classification result.

\begin{thm}
Let $\uC = (C_1,\ldots,C_d)$ be a row contraction on a Hilbert space $\CH_C$.
Then there is a one-to-one correspondence between
\begin{itemize}
\item[(a)]
unitary equivalence classes of subisometric liftings $\uE$ of $\uC$,
\item[(b)]
$\underline{L \otimes \eins}$-invariant subspaces $\CN$ of $\Gamma \otimes \CD_C$,
\item[(c)]
multi-analytic inner operators $M$ with symbols
$\Theta: \CD \rightarrow \Gamma \otimes \CD_C$ up to equivalence.
\end{itemize}
The correspondence is described above. In particular if $\uE$ is the lifting
then $\CD = \CD_E$, $M = M_{C,E}$ with symbol $\Theta = \Theta_{C,E}$ and $(b) \leftrightarrow (c)$ is Beurling's theorem.
\end{thm}

Theorem 1.6 shows that the characteristic function of a subisometric lifting characterizes the lifting up to unitary equivalence, justifying to call it characteristic.

\begin{proof}
$(b) \leftrightarrow (c)$ is Beurling's theorem, see \cite{Po89b}. We now show that the
correspondence $(a) \rightarrow (c)$ is well defined. Let $\uE$ on $\CH_E \supset \CH_C$
and $\uE^\prime$ on $\CH_{E^\prime} \supset \CH_C$ be two subisometric liftings of $\uC$
which are unitarily equivalent, i.e., there exists a unitary $u: \CH_E \rightarrow \CH_{E^\prime}$
such that $u |_{\CH_C} = \eins |_{\CH_C}$ and $E^\prime_i u = u E_i$ for $i=1,\ldots,d$.
Clearly unitarily equivalent row contractions have unitarily equivalent mids and we can
extend $u$ (in a trivial way) to a unitary $\hat{u}$ between the spaces $\hat{\CH}_E$
and $\hat{\CH}_{E^\prime}$ of the mids $\uV^E$ and $\uV^{E^\prime}$, so we have
\begin{align*}
\hat{u}: \hat{\CH}_E \rightarrow \hat{\CH}_{E^\prime} \;\mbox{unitary},
\quad \hat{u} |_{\CH_E} = u, \quad V^{E^\prime}_i \hat{u} = \hat{u} V^E_i \quad(i=1,\ldots,d)
\end{align*}
Because $\uE,\uE^\prime$ are subisometric we also have unitaries $W,W^\prime$ such that
\begin{align*}
W: \hat{\CH}_E \rightarrow \hat{\CH}_C, \quad V^C_i W = W V^E_i,
\quad W |_{\CH_C} = \eins |_{\CH_C}
\\
W^\prime: \hat{\CH}_{E^\prime} \rightarrow \hat{\CH}_C, \quad V^C_i W^\prime = W^\prime V^{E^\prime}_i, \quad W^\prime |_{\CH_C} = \eins |_{\CH_C} \nonumber
\end{align*}
If we now define
\[
u_C := W^\prime \hat{u} W^*: \hat{\CH}_C \rightarrow \hat{\CH}_C
\]
then it follows that $u_C$ commutes with the $V^C_i$ for
$i=1,\ldots,d$. To see that, ``chase" the following commuting diagram
\begin{align}\label{1.21}
\xymatrix@!{
 & \hat{\CH}_C \ar[rr]^{u_C} \ar'[d][dd]_{V^C_i}
   & & \hat{\CH}_C \ar[dd] \ar'[d][dd]_{V^C_i}
\\
 \hat{\CH}_E \ar[ur]^{W}
\ar[rr]^{\quad\quad\quad\quad\hat{u}} \ar'[d]^{V^E_i}[dd] \ar[dd]
 & & \hat{\CH}_{E'} \ar[ur]^{W^\prime} \ar'[d]^{V^{E^\prime}_i}[dd] \ar[dd] &
\\
 & \hat{\CH}_C \ar'[r]_{u_C}[rr]
   & & \hat{\CH}_C
\\
   \hat{\CH}_E \ar[rr]_{\hat{u}} \ar[ur]_W
 & & \hat{\CH}_{E'} \ar[ur]_{W^\prime}
}
\end{align}
Further, because $W, W^\prime$ and $\hat{u}$ all fix $\CH_C$ pointwise the same is true
for $u_C$, so we have also $u_C |_{\CH_C} = \eins |_{\CH_C}$. But by minimality of $\uV^C$
we know that $\hat{\CH}_C$ is the closed linear span of vectors of the form $V^C_\alpha h$
with $\alpha \in \tilde{\Lambda},\,h \in \CH_C$ and from
\[
u_C V^C_\alpha h = V^C_\alpha u_C h = V^C_\alpha h
\]
we infer that $u_C = \eins$.
Hence $W = (u_C)^* W^\prime \hat{u} = W^\prime \hat{u}$. Clearly $\hat{u}$
maps $e_0 \otimes \CD_E \subset \hat{\CH}_E$
onto $e_0 \otimes \CD_{E^\prime} \subset \hat{\CH}_{E^\prime}$,
so if we define the unitary $v := \hat{u} |_{\CD_E}: \CD_E \rightarrow \CD_{E^\prime}$
and use that $\Theta = W |_{\CD_E}$ and $\Theta^\prime = W^\prime |_{\CD_{E^\prime}}$ we see that
$\Theta = \Theta^\prime v$, i.e., the characteristic functions are equivalent.

Conversely suppose that we are given a multi-analytic inner operator with symbol
$\Theta: \CD \rightarrow \Gamma \otimes \CD_C$, as in (c). By $(b) \leftrightarrow (c)$ (Beurling's theorem) we have an invariant subspace $\CN$ which is associated to a subisometric lifting $\uE$ of $\uC$ and $\CD = \CD_E$, see the discussion preceding the theorem.
It remains to show that if $\Theta = \Theta^\prime v$ with a unitary $v: \CD_E \rightarrow \CD_{E^\prime}$ for two subisometric liftings $\uE$ and $\uE^\prime$ then $\uE$ and $\uE^\prime$
are unitarily equivalent. Let $W,W^\prime$ be the corresponding unitaries from the subisometric
lifting property. Then
\begin{align*}
W^\prime \CH_{E^\prime} \;=\;  \CH_C \oplus ( \Gamma \otimes \CD_C )
\;\ominus\; W^\prime ( \Gamma \otimes \CD_{E^\prime} ) \\
= \CH_C \oplus ( \Gamma \otimes \CD_C )
\;\ominus\; M_{C,E^\prime} \big( \Gamma \otimes v \CD_E \big) \\
= \CH_C \oplus ( \Gamma \otimes \CD_C )
\;\ominus\; M_{C,E} ( \Gamma \otimes \CD_E )
\;=\; W \CH_E,
\end{align*}
and we can define
\begin{align*}
U := (W^\prime)^* \, W |_{\CH_E}: \CH_E \rightarrow \CH_{E^\prime}.
\end{align*}

Because for $h \in \CH_C,$ $W h = h = W^\prime h$ we have $U h = h$.
In general for $h\in \CH_E$ and $i=1,\ldots,d$ (with $p_E, p_{E^\prime}$ orthogonal
projections onto $\CH_E, \CH_{E^\prime}$)
\begin{align*}
U E_i\, h = (W^\prime)^*\, W \,E_i \,h = (W^\prime)^* \,W \,p_E \,V^E_i \,h
= p_{E^\prime}\, (W^\prime)^*\, W\, V^E_i \,h
\\
= p_{E^\prime}\, (W^\prime)^*\, V^C_i\, W\, h
= p_{E^\prime}\, V^{E^\prime}_i\, (W^\prime)^* \,W\, h = E^\prime_i \,U\, h,
\end{align*}
i.e., $\uE$ and $\uE^\prime$ are unitarily equivalent.
\end{proof}

There is an interesting variant of the classification if we not only give $\uC$ but also $\uA$,
i.e., if we consider liftings of $\uC$ by $\uA$.

\begin{thm}
Let $\uA$ and $\uC$ be row contractions, $\uA$ $*$-stable.
There is a one-to-one correspondence between
\begin{itemize}
\item[(a)]
unitary equivalence classes of subisometric liftings of $\uC$ by $\uA$
\item[(b)]
equivalence classes of isometries $\gamma: \CD_{*,A} \rightarrow \CD_C$,
two isometries considered equivalent if they have the same range
\end{itemize}
\end{thm}

\begin{proof}
The details of this correspondence have already been discussed in connection with Proposition 1.2.
It is shown there how to construct an isometry $\gamma: \CD_{*,A} \rightarrow \CD_C$ if
a subisometric lifting of $\uC$ by $\uA$ is given, and conversely, how to use such an isometry
to find a subisometric lifting. The equivalence in (b) is chosen in such a way that two isometries
are equivalent if and only if the associated invariant subspaces are the same, compare
(1.9) and (1.11). Hence the result follows from Theorem 1.6.
\end{proof}

\begin{cor}
Let $\uA$ and $\uC$ be row contractions, $\uA$ $*$-stable.
A subisometric lifting of $\uC$ by $\uA$ exists
if and only if
\begin{align*}
dim \CD_{*,A} \leq dim \CD_C,
\end{align*}
where $dim$ stands for the cardinality of an orthonormal basis.
In the case $dim \CD_{*,A} = dim \CD_C$ (minimal subisometric dilation in the
terminology of \cite{DF84}) the lifting is unique up to unitary equivalence.
\end{cor}

\section{Coisometric Liftings}

The theory of subisometric liftings turns out to be especially relevant in the case of coisometric row contractions and coisometric liftings. We start with definitions and
elementary properties.

A row contraction $\uC$ on $\CH_1$ is called {\it coisometric} if $\uC \,\uC^* = \eins$, i.e.,
$\sum^d_{i=1} C_i C^*_i = \eins$. It is easy to check that a lifting $\uE$ on
${\cal H} = {\cal H}_C \oplus {\cal H}_A$ with block matrices
\begin{align*}
E_i =\left(
\begin{array}{cc}
C_i  & 0 \\
B_i  & A_i \\
\end{array}
\right)
\end{align*}
(for all $i=1,\ldots,d$) is coisometric if and only if $\underline{C}$ is coisometric and
\begin{align}\label{2.1}
\uB\, \uC^* = 0,\quad \mbox{i.e.},\; \sum^d_{i=1} B_i C^*_i = 0,
\\
\uA\, \uA^* + \uB\, \uB^* = \eins,\quad \mbox{i.e.},\;
\sum^d_{i=1} A_i A^*_i + \sum^d_{i=1} B_i B^*_i = \eins.
\end{align}

\begin{thm}
Let $\uA$ and $\uC$ be row contractions, $\uC$ coisometric.
Then there is a one-to-one correspondence between
\begin{itemize}
\item[(a)]
coisometric liftings $\uE$ of $\uC$ by $\uA$
\item[(b)]
isometries $\gamma: \CD_{*,A} \rightarrow \CD_C$
\end{itemize}

Explicitly, if
$E_i = \left(
\begin{array}{cc}
C_i  & 0 \\
B_i  & A_i \\
\end{array}
\right)$
for $i=1,\ldots,d$
provides a coisometric lifting $\uE$ of $\uC$ by $\uA$
then $\gamma: D_{*,A} h \mapsto \uB^* h \subset \CD_C$ (for $h \in \CH_A$)
is isometric.
\\

Conversely, if $\gamma: \CD_{*,A} \rightarrow \CD_C$ is isometric then with
$\uB^* := \gamma D_{*,A}$ we obtain a coisometric lifting $\uE$
by
$E_i = \left(
\begin{array}{cc}
C_i  & 0 \\
B_i  & A_i \\
\end{array}
\right)$
for $i=1,\ldots,d$.
\end{thm}

\begin{proof}
Because $\uC$ is coisometric, $D_C=\eins -\uC^* \uC$ is the orthogonal projection onto the kernel of $\uC$.

Let $\uE$ be a coisometric lifting of $\uC$ by $\uA$. Then from (2.1)
we have $\uC \uB^* = (\uB \uC^*)^* =0 $ and hence $range(\uB^*) \subset \CD_C$.

Further for $h \in \CH_A$, using (2.2)
\begin{align*}
\| D_{*,A} h \|^2 = \langle (\eins - \uA\uA^*)h, h \rangle =\langle \uB\uB^*h, h \rangle = \| \uB^*h\|^2
\end{align*}
So there exist an isometry $\gamma : \CD_{*,A} \to range(\uB^*) \subset \CD_C$ with
$\gamma D_{*,A} h = \uB^*h$ for all $h \in \CH_A$.

Conversely, let $\gamma: \CD_{*,A} \rightarrow \CD_C$ be an isometry and define
$\uB^* := \gamma D_{*,A}$. From $\uC |_{\CD_C} = 0$ we obtain $\uC\, \uB^* = 0$ or
$\uB \,\uC^*=0 $, which is (2.1). Further
\begin{align*}
\uB \,\uB^* = D_{*,A} \gamma^* \gamma D_{*,A} = D^2_{*,A} = \eins - \uA\, \uA^*,
\end{align*}
hence $\uA\, \uA^* + \uB\, \uB^* = \eins$, which is (2.2). Hence with
$E_i = \left(
\begin{array}{cc}
C_i  & 0 \\
B_i  & A_i \\
\end{array}
\right),$
for $i=1,\ldots,d,$
we obtain a coisometric lifting $\uE$ of $\uC$ by $\uA$.

Finally if $\gamma, \gamma^\prime: \CD_{*,A} \rightarrow \CD_C$ are two isometries and
$\gamma \not= \gamma^\prime$ then $\uB^* \not= (\uB^\prime)^*$ for $\uB^*,(\uB^\prime)^*$
defined by $\gamma, \gamma^\prime$ as above. Hence the correspondence is one-to-one.
\end{proof}

\begin{cor}
Let $\uA$ and $\uC$ be row contractions, $\uC$ coisometric.
A coisometric lifting $\uE$ of $\uC$ by $\uA$ exists
if and only if
\begin{align*}
dim \CD_{*,A} \leq dim \CD_C,
\end{align*}
where $dim$ stands for the cardinality of an orthonormal basis..
\end{cor}

Theorem 2.1 gives a kind of free parametrization of the coisometric liftings.
Let us consider an elementary example.
\begin{align} \label{2.3}
\uc = (c_1,\ldots,c_d) \in \C^d, \quad \|\uc\|^2 = \sum^d_{i=1} |c_i|^2 = 1
\quad {\mbox (unit\, sphere)}
\\
\ua = (a_1,\ldots,a_d) \in \C^d, \quad \|\ua\|^2 = \sum^d_{i=1} |a_i|^2 \leq 1
\qquad {\mbox (unit\, ball)} \nonumber
\end{align}
Then we get a left lower corner $\ub = (b_1,\ldots,b_d)$ for a coisometric lifting
if $\langle \ub,\uc \rangle = 0$ and $\|\ua\|^2 + \|\ub\|^2 = 1$, according to (2.1) and (2.2).
Obviously the set of solutions for $\ub$ is the (complex) sphere with radius
$r = \sqrt{1-\|\ua\|^2}$ in the subspace orthogonal to $\uc$. If $\|\ua\|=1$ the
solution is unique.
We can check that the parametrization using isometries $\gamma: \CD_{*,A} \rightarrow \CD_C$
as in Theorem 2.1 yields the same result.
\\

Theorem 2.1 and Corollary 2.2 are even true if $\uA$ is not $*$-stable. If $\uA$ is $*$-stable
then we should compare these results with those in Section 1. Note in particular that
the formula $\uB^* = \gamma D_{*,A}$ in Theorem 2.1 and the formula
$\uB^* = D^*_C \,\gamma\, D_{*,A}$ (1.10) are compatible because, as noted above,
for $\uC$ coisometric the operator $D^*_C$ is nothing but the embedding of $\CD_C$ into
$\oplus^d_{i=1} \CH_C$ which is implicit in the formulation chosen in Theorem 2.1.
Further comparison yields the following result
which shows that subisometric liftings occur very naturally in the coisometric setting.

\begin{prop}
Let $\uC$ be a coisometric row contraction. A lifting of $\uC$ is a coisometric lifting by a $*$-stable $\uA$ if and only if it is subisometric.
\end{prop}

\begin{proof}
Using Theorem 2.1 we can replace the condition ``coisometric" for
the lifting by the existence of an isometry $\gamma: \CD_{*,A}
\rightarrow \CD_C$ such that $\uB^* = \gamma D_{*,A} = D^*_C \gamma
D_{*,A}$. Now Proposition 2.3 is a direct consequence of Proposition
1.2.
\end{proof}

In particular, for coisometric liftings by a $*$-stable $\uA$ there
exists an associated invariant subspace and a characteristic
function. In the special case $dim\, \CH_C = 1$ this characteristic
function was introduced in \cite{DG07} under the name ``extended
characteristic function". For general $\CH_C$, in view of Theorem
1.6, it is better to call it the characteristic function of the
lifting (with $\uC$ given), as we have done in Definition 1.5.

\section{Characteristic Functions of Reduced Liftings}
In this section we generalize the theory of characteristic functions for sub\-isometric liftings from Section 1 and establish a setting that also includes the setting of Section 2.

Let $\uC$ be a row contraction on $\CH_C$ and $\uE$ on $\CH_E = \CH_C \oplus \CH_A$ be a (contractive) lifting so that for all $i=1,\ldots,d$
\begin{align*}
E_i =\left(
\begin{array}{cc}
C_i  & 0 \\
B_i  & A_i \\
\end{array}
\right)
\end{align*}
Then as in \eqref{1.3} we have a mid $\uV^E$ on $\CH_E \oplus \big( \Gamma \otimes \CD_E \big)$. Clearly $\uV^E$ is an isometric lifting of $\uC$, so the space of the mid $\uV^C$
can be embedded as a subspace reducing the $V^E_i$. Let us encode this by introducing the restriction $\uY$ on the orthogonal complement $\CK$ and a unitary $W$ by
\begin{align}\label{3.1}
W: \CH_E \oplus \big( \Gamma \otimes \CD_E \big)
\rightarrow
\CH_C \oplus \big( \Gamma \otimes \CD_C \big) \oplus \CK \\
\tilde{V}^E_i W = W V^E_i, \quad W |_{\CH_C} = \eins |_{\CH_C} \quad
\mbox{with} \quad \tilde{\uV}^E = \uV^C \oplus \uY \nonumber
\end{align}
By omitting $\CH_C$ we also have a unitary (also denoted by $W$)
\begin{align}\label{3.2}
W: \CH_A \oplus \big( \Gamma \otimes \CD_E \big)
\rightarrow
\big( \Gamma \otimes \CD_C \big) \oplus \CK
\end{align}
and an isometric embedding $\CH_{A_*} := W \CH_A \subset \big( \Gamma \otimes\CD_C \big) \oplus \CK$.
Further we obtain
\begin{align}\label{3.3}
\uB^* = p_{\CH_C} (\uV^E)^* |_{\CH_A} = p_{\CH_C} \big[(\uV^C)^* \oplus \uY^* \big] W |_{\CH_A} = D^*_C\, p_{e_0 \otimes \CD_C} W |_{\CH_A}
\end{align}
where we used formula \eqref{1.3} for $\uV^C$.
\\

To proceed we need a few facts about the mid $\uV^A$ on $\tilde{\CH}_A$ of $\uA$. We write its Wold decomposition as
\begin{align}\label{3.4}
\tilde{\CH}_A = \big( \Gamma \otimes \CD_{*,A} \big) \oplus \CR_A \\
V^A_i = (L_i \otimes \eins) \oplus R^A_i, \quad i=1,\ldots,d, \nonumber
\end{align}
where $\CR_A$ and $\uR^A$ stand for the residual part (cf. \cite{Po89a}).
The embedding of $\CH_A$ into $\tilde{\CH}_A$ can be written as
\begin{align}\label{3.5}
\CH_A \ni h \mapsto \big( \sum_{\alpha \in \tilde{\Lambda}}
e_\alpha \otimes D_{*,A} A^*_\alpha h \big) \oplus h_{\CR}
\end{align}
Here $h_{\CR}$ belongs to the residual part $\CR_A$.
Compare \cite{BDZ06} for a derivation of this decomposition via Stinespring's theorem. In fact,
it is not difficult to check that a formula like \eqref{3.5} always reproduces the Wold decomposition above, compare also \cite{FF90} for similar arguments for $d=1$.
Note that the residual part vanishes if and only if $\uA$ is $*$-stable, so in this case
we are back in the setting of Section 1.

Further we need the decomposition $\CH_A = \CH^1_A \oplus \CH^0_A$ with $\CH^1_A$
the largest subspace invariant for the $A^*_i$ and such that the restriction of $\uA^*$ is isometric, i.e.,
\begin{align}\label{3.6}
\CH^1_A := \{ h \in \CH_A: \sum_{|\alpha|=n} \| A^*_\alpha h \|^2 = \|h\|^2
\;\mbox{for all}\; n \in \N \}
\end{align}
Then it is easy to check that $\CH^1_A = \CH_A \cap \CR_A$ (cf. \cite{Po89a}, Proposition 2.9), but the position of
$\CH^0_A$ may be complicated with respect to the decomposition \eqref{3.5} because $\uA$
restricted to $\CH^0_A$ may not be $*$-stable and in this case $\CH^0_A$ is not contained in
$\Gamma \otimes \CD_{*,A}$. In fact, if $0 \not= h \in \CH^0_A$ we only have
\begin{align}\label{3.7}
\sum_{|\alpha|=n} \| A^*_\alpha h \|^2 < \|h\|^2
\;\mbox{for some}\; n \in \N
\end{align}
which (by definition) means that $\uA |_{\CH^0_A}$ is {\em
completely non-coisometric (c.n.c.)}, cf. \cite{Po89a}.

Now we look at $\uA$ and its mid $\uV^A$ embedded into the larger
structure obtained from the lifting $\uE$. Clearly $\uV^E$
restricted to $\CH_A \oplus \big( \Gamma \otimes \CD_E \big)$ is an
isometric dilation of $\uA$, so $\CH_A \oplus \big( \Gamma \otimes
\CD_E \big)$ contains $\tilde{\CH}_A$ as a $V^E_i$-reducing subspace
($i=1,\ldots,d$) which we still denote by $\tilde{\CH}_A$. Using
\eqref{3.2} we see that $\big( \Gamma \otimes \CD_C \big) \oplus
\CK$ contains the $(L_i \otimes \eins) \oplus Y_i$-reducing subspace
$W \tilde{\CH}_A$ and the restriction of $(\underline{L \otimes
\eins}) \oplus \uY$ is a mid of $\uA$ (transferred to $W \CH_A$).
Denoting the restriction of $W$ to $\tilde{\CH}_A$ also by $W$
we have (for $i=1,\ldots,d$)
\begin{align}\label{3.8}
W \big[ (L_i \otimes \eins) \oplus R^A_i \big] = W V^A_i
= \big[ (L_i \otimes \eins) \oplus Y_i \big] W
\end{align}
Where is $\CH_{A_*} = W \CH_A$? Clearly
\begin{align}\label{3.9}
W \CH^1_A = W (\CH_A \cap \CR_A) \subset W \CR_A \subset \CK,
\end{align}
where the last inclusion follows from \eqref{3.8} and the fact that $\underline{L \otimes \eins}$
is $*$-stable. The position of $W \CH^0_A$ may be more complicated.

To organize the relevant data we use \eqref{3.4} together with the embedding of
$\tilde{\CH}_A$ into $\CH_A \oplus \big( \Gamma \otimes \CD_E \big)$
and \eqref{3.2} to define
\begin{align}\label{3.10}
M: \Gamma \otimes \CD_{*,A} \rightarrow \Gamma \otimes \CD_C, \\
M = P_{\Gamma \otimes \CD_C} W |_{\Gamma \otimes \CD_{*,A}} \nonumber
\end{align}
which is a multi-analytic operator. Then for $h \in \CH_A$
\begin{align*}
P_{e_0 \otimes \CD_C} W h = P_{e_0 \otimes \CD_C} M P_{\Gamma \otimes \CD_{*,A}} h
= P_{e_0 \otimes \CD_C} M P_{e_0 \otimes \CD_{*,A}} h
\end{align*}
where for the first equality we used \eqref{3.9} and the second then follows from the fact that
$M$ is multi-analytic. But $P_{e_0 \otimes \CD_{*,A}} h = e_0 \otimes D_{*,A} h$ by \eqref{3.5}
and we conclude that $P_{e_0 \otimes \CD_C} W |_{\CH_A}: \CH_A \rightarrow \CD_C$ factors through
$\CD_{*,A}$ in the sense that there exists a contraction
$\gamma := P_{e_0 \otimes \CD_C} M |_{e_0 \otimes \CD_{*,A}}:
\CD_{*,A} \rightarrow \CD_C$ such that
\begin{align}\label{3.11}
P_{e_0 \otimes \CD_C} W |_{\CH_A} = \gamma D_{*,A}
\end{align}
In fact, $\gamma$ is nothing but the the $0$-th Fourier coefficient of $M$ in the sense
of \cite{Po03}. Combined with \eqref{3.3} we obtain
\begin{align}\label{3.12}
\uB^* = D^*_C\, \gamma\, D_{*,A}: \CH_A \rightarrow \bigoplus^d_{i=1} \CH_C
\end{align}
This is one half of the following analogue for row contractions of Lemma 2.1 in Chap.IV of
\cite{FF90}, which already has been discussed in the introduction, see in particular \eqref{0.3}.
\begin{prop}
$\uE = (E_1,\ldots,E_d)$ on $\CH_E = \CH_C \oplus \CH_A$ with block matrices
\begin{align*}
E_i =\left(
\begin{array}{cc}
C_i  & 0 \\
B_i  & A_i \\
\end{array}
\right)
\end{align*}
(for $i=1,\ldots,d$) is a row contraction if and only if $\uC$ and $\uA$ are row contractions
and there exists a contraction $\gamma: \CD_{*,A} \rightarrow \CD_C$ such that
\eqref{3.12} holds.
\end{prop}

\begin{proof}
Clearly if $\uE$ is a row contraction then $\uC$ and $\uA$ are row contractions. Above we have already given a (dilation) proof that if $\uE$ is contractive then $\uB$ satisfies \eqref{3.12} for a suitable contraction $\gamma$. To prove the converse,
let $\gamma: \CD_{*,A} \rightarrow \CD_C$ be a contraction and $\uB^*$ given as in \eqref{3.12}.
Then for $x \in \CH_C,\; y \in \CH_A$
\begin{align*}
| \langle x, \uC\, \uB^* y \rangle |^2 = | \langle x, \uC D^*_C \gamma D_{*,A}\, y \rangle |^2
= | \langle D_C \uC^* x, \gamma D_{*,A}\, y \rangle |^2
\end{align*}
\begin{align*}
\leq \| D_C \uC^* x \|^2 \| \gamma D_{*,A} y \|^2
\leq \langle x, (\eins - \uC \uC^*) x \rangle \; \langle y, (\eins - \uA \uA^*) y \rangle
\end{align*}
which implies (see for example Exercise 3.2 in \cite{Pau03}) that
\begin{align*}
0 \leq
\left(
\begin{array}{cc}
\eins - \uC \uC^* & -\uC \uB^*\\
-\uB \uC^*  & \eins - \uA \uA^* \\
\end{array}
\right)
= \eins - \uE \uE^*
\end{align*}
hence $\uE$ is a row contraction.
\end{proof}

Let us go back to the lifting $\uE$ of $\uC$ by $\uA$. The following definition is useful to analyze further the position of $W \CH_A$.
\begin{defn}
$\gamma: \CD_{*,A} \rightarrow \CD_C$ is called resolving if for all $h \in \CH_A$ we have
\begin{align*}
\big( \gamma D_{*,A} A^*_\alpha h =0 \;\mbox{for all}\; \alpha \in \tilde{\Lambda} \big)
\Rightarrow \big( D_{*,A} A^*_\alpha h =0 \;\mbox{for all}\; \alpha \in \tilde{\Lambda} \big)
\end{align*}
\end{defn}
Clearly if $\gamma: \CD_{*,A} \rightarrow \CD_C$ is injective then it is resolving.
Note that $D_{*,A} A^*_\alpha h =0 \;\mbox{for all}\; \alpha \in \tilde{\Lambda}$ if and only if
$h \in \CH^1_A$, and so
the intuitive meaning of `resolving' is that `\,looking at $\CH_A$ through $\gamma$\,' still allows to detect whether $h \in \CH_A$ is in $\CH^1_A$ or not.
More precisely, $\gamma$ is resolving if and only if for all $h \in \CH^0_A = \CH_A \ominus \CH^1_A$
there exists $\alpha \in \tilde{\Lambda}$ such that $\gamma D_{*,A} A^*_\alpha h \not=0$.
In particular if $\uA$ is c.n.c., i.e. $\CH^1_A = \{0\}$, then $\gamma$ is resolving if and only if for all $0 \not= h \in \CH_A$ there exists $\alpha \in \tilde{\Lambda}$ such that $\gamma D_{*,A} A^*_\alpha h \not=0$.
\begin{lem}
The following assertions are equivalent
\begin{itemize}
\item[(a)]
$\gamma$ is resolving.
\item[(b)]
$W \CH_A \cap \CK \subset W \CH^1_A$
\item[(c)]
$W \CH_A \cap \CK \;=\; W \CH^1_A$
\item[(d)]
$\big( \Gamma \otimes \CD_C \big) \bigvee W \big( \Gamma \otimes \CD_E \big)
\; = \; \big( \Gamma \otimes \CD_C \big) \oplus (\CK \ominus W \CH^1_A)$
\end{itemize}
\end{lem}
\begin{proof}
(b) says that for $h \in \CH_A \setminus \CH^1_A$ the embedded
$Wh$ is not in $\CK$, so not orthogonal to $\Gamma \otimes \CD_C$,
equivalently, there exists $\alpha \in \tilde{\Lambda}$ such that
\begin{align*}
0 \not= P_{e_0 \otimes \CD_C} \big[ (L^*_\alpha \otimes \eins) \oplus Y^*_\alpha \big] W h
= P_{e_0 \otimes \CD_C} W (V^A_\alpha)^* h = \gamma D_{*,A} A^*_\alpha h
\end{align*}
(where we used the embedding of the mid of $\uA$ and in particular \eqref{3.11}).
By comparison with the comments following Definition 3.2 we conclude that (a) and (b) are equivalent. We noted in \eqref{3.9} that always
$W \CH^1_A \subset \CK$, so (b) and (c) are equivalent.

To get the equivalence of (c) and (d) note that $x \in \big( \Gamma \otimes \CD_C \big) \oplus \CK$ is orthogonal to $\big( \Gamma \otimes \CD_C \big)$ and to
$W \big( \Gamma \otimes \CD_E \big)$ if and only if $x \in \CK$ and $x \in W \CH_A$
(compare \eqref{3.2}).
Hence the orthogonal complement of $\big( \Gamma \otimes \CD_C \big) \bigvee W \big( \Gamma \otimes \CD_E \big)$ in $\big( \Gamma \otimes \CD_C \big) \oplus \CK$ is in fact $W \CH_A \cap \CK$.
\end{proof}
\begin{defn}
A lifting $\uE$ of $\uC$ by $\uA$ is called reduced if
$\uA$ is c.n.c. (i.e., $\CH^1_A = \{0\}$, see \eqref{3.7}) and $\gamma$ is resolving.
\end{defn}
We have already seen two important classes of reduced liftings.
\begin{itemize}
\item[1)]
Subisometric liftings. Here $\uA$ is $*$-stable and $\gamma$ is isometric, see Proposition 1.2.
\item[2)]
Coisometric liftings by $\uA$ c.n.c. Here $\gamma$ is isometric by Theorem 2.1.
\end{itemize}
Note that by Proposition 2.3 the coisometric liftings by $*$-stable $\uA$ are exactly the intersection of cases 1) and 2).
\begin{lem}
The following assertions are equivalent
\begin{itemize}
\item[(a)]
$\uE$ is reduced.
\item[(b)]
$\{ h \in \CH_A: \gamma D_{*,A} A^*_\alpha h = 0 \;\mbox{for all}\; \alpha \in \tilde{\Lambda}\} \;=\; \{0\}$
\item[(c)]
$W \CH_A \cap \CK \;=\; \{0\}$
\end{itemize}
\end{lem}
\begin{proof}
If $\gamma$ is resolving then (by definition) the space given in (b) is contained in $\CH^1_A$.
Hence (a) implies (b). Also, from (b) we first conclude that
$\CH^1_A = \{ h \in \CH_A: D_{*,A} A^*_\alpha h = 0 \;\mbox{for all}\; \alpha \in \tilde{\Lambda}\} =\; \{0\}$ and then that $\gamma$ is resolving, so (b) implies (a).
If we have (c) then by Lemma 3.3(b) $\gamma$ is resolving and then by Lemma 3.3(c) $\uA$
is c.n.c., so we have (a). Given (a), Lemma 3.3(c) implies (c).
\end{proof}

If $\gamma D_{*,A} A^*_\alpha h = 0 \;\mbox{for all}\; \alpha \in \tilde{\Lambda}$
then by \eqref{3.12} we conclude that $A^*_\alpha h \in ker \uB^* = (range \uB)^\perp$.
Hence vectors in the space
$\{ h \in \CH_A: \gamma D_{*,A} A^*_\alpha h = 0 \;\mbox{for all}\; \alpha \in \tilde{\Lambda}\}$
do not contribute in any way
to the interaction between $\CH_A$ and $\CH_C$ via $\uB^*$, and it is no great loss to concentrate on liftings where this space has been removed. By Lemma 3.5(b), in doing this we obtain
exactly the reduced liftings. This also explains our terminology.
\\

For reduced liftings we can successfully develop a theory of characteristic functions.
\begin{defn}
Let $\uE$ be a reduced lifting of $\uC$ by $\uA$. We call the multi-analytic operator
\begin{align}\label{3.13}
M_{C,E}: \Gamma \otimes \CD_E \rightarrow \Gamma \otimes \CD_C, \\
M_{C,E} = P_{\Gamma \otimes \CD_C} W |_{\Gamma \otimes \CD_E} \nonumber
\end{align}
(or its symbol $\Theta_{C,E}: \CD_E \rightarrow \Gamma \otimes \CD_C$) the characteristic function of the lifting $\uE$.
\end{defn}


Using the characteristic function we can develop a theory of functional models for reduced liftings. The idea is similar as in the case of characteristic functions for c.n.c. row contractions,
see \cite{Po89b}.

Let $\uE$ be a reduced lifting of $\uC$ by $\uA$. From $\uA$ c.n.c. we obtain $\CH^1_A = \{0\}$
and then Lemma 3.3 gives
\begin{align} \label{3.14}
\big( \Gamma \otimes \CD_C \big) \bigvee W \big( \Gamma \otimes \CD_E \big)
\; = \; \big( \Gamma \otimes \CD_C \big) \oplus \CK
\end{align}
With the definition
\begin{align} \label{3.15}
\Delta_{C,E} := (\eins - M^*_{C,E} M_{C,E})^{\frac{1}{2}}:
\Gamma \otimes \CD_E \rightarrow \Gamma \otimes \CD_E
\end{align}
we obtain for $x \in \Gamma \otimes \CD_E$
\begin{align} \label{3.16}
\|P_{\CK} W x \|^2 = \| (\eins - P_{\Gamma \otimes \CD_C}) W x \|^2
= \|x\|^2 - \| P_{\Gamma \otimes \CD_C} W x \|^2
\\
= \|x\|^2 - \| M_{C,E} x \|^2 = \| \Delta_{C,E} x \|^2 \nonumber
\end{align}
This means that we can isometrically identify $\CK$ with
$\overline{\Delta_{C,E}(\Gamma \otimes \CD_E)}$ and with
this identification we have
\begin{align} \label{3.17}
W \CH_A = \big[ (\Gamma \otimes \CD_C) \oplus \CK  \big]
\ominus W (\Gamma \otimes \CD_E)
\end{align}
\begin{align*}
= \big[ (\Gamma \otimes \CD_C) \oplus \overline{\Delta_{C,E}(\Gamma
\otimes \CD_E)} \big] \ominus \{ M_{C,E}\, x \oplus \Delta_{C,E}\,
x: x \in \Gamma \otimes \CD_C \} \nonumber
\end{align*}
which is a kind of functional model.
\begin{thm}
Let $\uC$ be a row contraction. Reduced liftings $\uE$ and $\uE'$ of $\uC$ are unitarily equivalent if and only if their characteristic functions $M_{C,E}$ and $M_{C,E'}$ are equivalent.
\end{thm}
Recall that $M_{C,E}$ and $M_{C,E'}$ are equivalent if there exists a unitary $v: \CD_E \rightarrow \CD_{E'}$ such that their symbols satisfy $\Theta_{C,E} = \Theta_{C,E'}\, v$.
Compared with the analogous result for subisometric liftings contained in Theorem 1.6 the modifications necessary to prove Theorem 3.7 are technical and straightforward, so we omit the
proof. The important thing to recognize is that, if a lifting $\uE$ is reduced, we have the
functional model \eqref{3.17} for it which is built only from $\uC$ and from the characteristic
function $M_{C,E}$.
\\

Conversely, if $\uC$ on $\CH_C$ is a row contraction and
\begin{align*}
\tilde{M}: \Gamma \otimes \CD \rightarrow \Gamma \otimes \CD_C
\end{align*}
is an arbitrary contractive multi-analytic function (where $\CD$ is any Hilbert space),
then we can define
\begin{align*}
\Delta :(\eins - \tilde{M}^* \tilde{M})^{\frac{1}{2}}: \;
\Gamma \otimes \CD \rightarrow \Gamma \otimes \CD
\end{align*}
\begin{align*}
\tilde{\CH} := \CH_C \oplus (\Gamma \otimes \CD_C) \oplus
\overline{\Delta (\Gamma \otimes \CD)}
\end{align*}
\begin{align*}
\tilde{W}: \Gamma \otimes \CD \rightarrow
(\Gamma \otimes \CD_C) \oplus \overline{\Delta (\Gamma \otimes \CD)},\; x \mapsto
\tilde{M}\, x \oplus \Delta \, x
\end{align*}
$\tilde{W}$ is isometric and by introducing a copy $\CH_A$ of the orthogonal complement
of $\tilde{W}(\Gamma \otimes \CD)$ we can extend $\tilde{W}$ to a unitary
\begin{align*}
\tilde{W}: \CH_A \oplus (\Gamma \otimes \CD) \rightarrow
(\Gamma \otimes \CD_C) \oplus \overline{\Delta (\Gamma \otimes \CD)}
\end{align*}
Let $\tilde{\uV} = (\tilde{V}_1,\ldots,\tilde{V}_d)$ be defined on $\tilde{\CH}$ by $\tilde{V}_i := V^C_i \oplus Y_i$ (for
$i=1,\ldots,d$), where $\uV^C$ is the mid of $\uC$ on $\CH_C \oplus (\Gamma \otimes \CD_C)$
\eqref{1.3} and $Y_i$ is given by
\begin{align*}
Y_i \Delta \,x := \Delta (L_i \otimes \eins) x \quad (\mbox{where} \; x \in \Gamma \otimes \CD)
\end{align*}
It is not difficult to check that $\uY$ (and hence also $\tilde{\uV}$) is a row contraction
consisting of isometries with orthogonal ranges (i.e., a row isometry). Further
\begin{align*}
\tilde{W}(\Gamma \otimes \CD) = \{ \tilde{M}\, x \oplus \Delta\,x,\; x \in \Gamma \otimes \CD \}
\end{align*}
is invariant for the $\tilde{V}_i$. With
$E^*_i := \tilde{V}^*_i |_{\CH_C \oplus \tilde{W} \CH_A}, \;
A^*_i := \tilde{V}^*_i |_{\tilde{W} \CH_A}$ for $i=1,\ldots,d$
we obtain a contractive lifting
$\uE$ of $\uC$ by $\uA$ which we may call {\it the lifting associated to the multi-analytic
function} $\tilde{M}$. The following result gives another justification for considering
reduced liftings.
\begin{prop}
The contractive lifting $\uE$ associated to a row contraction $\uC$ and a contractive multi-analytic
function $\tilde{M}: \Gamma \otimes \CD \rightarrow \Gamma \otimes \CD_C$
(where $\CD$ is any Hilbert space) is reduced.
\end{prop}
\begin{proof}
By Lemma 3.5 it is enough to show that any vector $y \in \tilde{W} \CH_A$ which is
orthogonal to $\Gamma \otimes \CD_C$ is the zero vector. But $y \in \tilde{W} \CH_A$
means that $y$ is orthogonal to $\tilde{M}\, x \oplus \Delta\,x$ for all
$x \in \Gamma \otimes \CD$ and $y$ orthogonal to $\Gamma \otimes \CD_C$
means that $y \in 0 \oplus \overline{\Delta (\Gamma \otimes \CD)}$. Hence indeed $y=0$.
\end{proof}

Proposition 3.8 shows that the theory of characteristic functions cannot be extended beyond
reduced liftings. Note that $\tilde{M}$ is not necessarily the characteristic function of
the associated lifting $\uE$ and we used $\,\tilde{}\,$ to indicate this. It is an
interesting question which intrinsic properties of $\tilde{M}$ guarantee that it is
the characteristic function. We leave this as an open problem.

\section{Properties of the Characteristic Function}

First we shall compute an explicit expression for the characteristic function of a reduced lifting. We continue to use the notation of the previous section and consider
a reduced lifting $\uE$ on
$\CH_E = \CH_C \oplus \CH_A$ of $\uC$ on $\CH_C$ by $\uA$ on $\CH_A$.
As in \eqref{3.8} the row isometry $(\underline{L \otimes \eins}) \oplus \uY$
on $(\Gamma \otimes \CD_C) \oplus \CK$ restricts to a mid of $\uA$ (transferred
to $W \CH_A$). So we have for all $\alpha \in \tilde{\Lambda}$ and $h \in \CH_A$
\begin{align} \label{4.1}
\big[ (L^*_\alpha \otimes \eins) \oplus Y^*_\alpha \big] W h = W \, A^*_\alpha h
\end{align}
Using \eqref{3.11} we infer that
\begin{align} \label{4.2}
\gamma D_{*,A} A^*_\alpha h
= P_{e_0 \otimes \CD_C} W \, A^*_\alpha h
= P_{e_0 \otimes \CD_C} \big[ (L^*_\alpha \otimes \eins) \oplus Y^*_\alpha \big] W h
= P_{e_\alpha \otimes \CD_C} W h
\end{align}
which yields a Poisson kernel type formula, compare \eqref{1.8}:
\begin{align} \label{4.3}
P_{\Gamma \otimes \CD_C} W h
= \sum_{\alpha \in \tilde{\Lambda}} e_\alpha \otimes \gamma D_{*,A} A^*_\alpha h
\end{align}
To compute the symbol $\Theta_{C,E}$ of the characteristic function we define
$d^i_h := (V^E_i - E_i) h = e_0 \otimes (D_E)_i h$ and use the identification of $\CD_E$ with the closed linear span of all $d^i_h$ with $i=1,\ldots,d$ and $h \in \CH_E$, see \eqref{1.3}.
Then, using \eqref{3.1} and the
Definition 3.6 of $\Theta_{C,E}$, we obtain
\begin{align} \label{4.4}
\Theta_{C,E} d^i_h = P_{\Gamma \otimes \CD_C} W (V^E_i - E_i) h
= P_{\Gamma \otimes \CD_C} V^C_i P_{\CH_C \oplus (\Gamma \otimes \CD_C)} W h - P_{\Gamma \otimes \CD_C} W E_i h
\end{align}
We distinguish two cases.

\vspace{0.2cm}

{\bf Case I:} $h \in \CH_C.$
\begin{align*}
P_{\Gamma \otimes \CD_C} V^C_i P_{\CH_C \oplus (\Gamma \otimes \CD_C)} W h
= P_{\Gamma \otimes \CD_C} V^C_i h
= [e_0 \otimes (D_C)_i h ] \quad\quad {\mbox{by} \;\eqref{1.3}}
\end{align*}
\begin{align*}
P_{\Gamma \otimes \CD_C} W E_i h
= P_{\Gamma \otimes \CD_C} W (C_i h \oplus B_i h)
= \sum_\alpha e_\alpha \otimes \gamma \DsA A^*_\alpha B_i h \quad\quad {\mbox{by} \;\eqref{4.3}}
\end{align*}
and thus
\begin{align} \label{4.5}
\Theta_{C,E}  d^i_h = e_0 \otimes
[(D_C)_i h - \gamma \DsA B_i h] - \sum_{|\alpha|\geq 1} e_\alpha \otimes \gamma \DsA
A^*_\alpha B_i h
\end{align}

{\bf Case II:} $h \in\, \CH_A$.
\begin{align*}
P_{\Gamma \otimes \CD_C} V^C_i P_{\CH_C \oplus (\Gamma \otimes \CD_C)} W h
= V^C_i P_{\Gamma \otimes \CD_C} W h
\end{align*}
\begin{align*}
= (L_i \otimes \eins) P_{\Gamma \otimes \CD_C} W h
= \sum_\alpha e_i \otimes e_\alpha \otimes \gamma \DsA A^*_\alpha h
\end{align*}
\begin{align*}
P_{\Gamma \otimes \CD_C} W E_i h = P_{\Gamma \otimes \CD_C} W A_i h
= \sum_\beta e_\beta \otimes \gamma \DsA A^*_\beta A_i h
\end{align*}
Note that for $h \in\, \CH_A$ we have $(D_A)_i h = (D_E)_i h$ (which we identify with $d^i_h$)
because $\uE$ is an extension of $\uA$.
With $P_j$ the orthogonal projection onto the $j-$th component we obtain
\begin{align} \label{4.6}
\Theta_{C,E}\, d^i_h = - e_0 \otimes \gamma \DsA A_i h
+ \sum^d_{j=1} e_j \otimes \sum_\alpha e_\alpha \otimes \gamma \DsA A^*_\alpha
(\delta_{ji} \eins -A^*_j A_i) h
\nonumber
\\
= - e_0 \otimes \gamma \uA (D_A)_i h
+ \sum^d_{j=1} e_j \otimes \sum_\alpha e_\alpha \otimes \gamma \DsA A^*_\alpha
P_j D_A (D_A)_i h \nonumber
\\
= - e_0 \otimes \gamma \sum^d_{j=1} A_j P_j d^i_h
+ \sum^d_{j=1} e_j \otimes \sum_\alpha e_\alpha \otimes \gamma \DsA A^*_\alpha P_j D_A d^i_h
\end{align}
We note that if $\gamma$ is omitted from \eqref{4.6} then we obtain exactly Popescu's definition
of the characteristic function of the (c.n.c.) row contraction $\uA$ as given in \cite{Po89b}.
Hence Case II is essentially
the characteristic function of $\uA$ , contractively embedded by $\gamma$. In a special case this has been observed in \cite{DG07} and, because this special case was subisometric and hence $\gamma$
isometric, $\Theta$ was called an extended characteristic function. \eqref{4.6} generalizes this idea.
\\

Let us now illustrate how the characteristic function factorizes for
iterated liftings. Assume that $\tilde{\uE}$ on $\CH_{\tilde{E}}$ is
a two step lifting of the row contraction  $\uC$ on $\CH_C$, i.e.,
$\uE$ on $\CH_E$ with $E_i= \left(\begin{array}{cc}
C_i & 0   \\
B_i & A_i
\end{array}
\right)$ (for $i=1,\ldots,d$) is a contractive lifting of $\uC$ on
$\CH_C$ by $\uA$ on $\CH_A$ (as before) and $\tilde{\uE}$ on
$\CH_{\tilde{E}}$ with $\tilde{E}_i= \left(\begin{array}{cc}
E_i & 0   \\
*   & \tilde{A}_i
\end{array}
\right)$
(for $i=1,\ldots,d$) is a contractive lifting of $\uE$ on $\CH_E$ by $\tilde{\uA}$ on $\CH_{\tilde{A}}$. Then
$\CH_{\tilde{E}} = \CH_E \oplus \CH_{\tilde{A}} = \CH_C \oplus \CH_A \oplus \CH_{\tilde{A}}$
and with respect to this decomposition
\begin{align} \label{4.7}
\tilde{E}_i=\left(\begin{array}{ccc}
C_i & 0   & 0 \\
*   & A_i & 0 \\
*   & *   & \tilde{A}_i
\end{array}
\right)
\end{align}
`$*$' stands for entries which we do not need to name explicitly.

\begin{thm}
If the liftings $\uE$ of $\uC$ and $\tilde{\uE}$ of $\uE$ are reduced
then also the lifting $\tilde{\uE}$ of $\uC$ is reduced, and the
characteristic functions factorize:
\begin{align} \label{4.8}
M_{C,\tilde{E}}=M_{C,E} \; M_{E,\tilde{E}}.
\end{align}
\end{thm}

\begin{proof}
As in \eqref{3.1} we obtain the following unitaries from the given liftings:
\begin{align*}
W: \CH_E \oplus (\Gamma \otimes \CD_E) \rightarrow \CH_C \oplus (\Gamma \otimes \CD_C) \oplus \CK
\end{align*}
\begin{align*}
\tilde{W}: \CH_{\tilde{E}} \oplus (\Gamma \otimes \CD_{\tilde{E}}) \rightarrow \CH_E \oplus (\Gamma \otimes \CD_E) \oplus \tilde{\CK}
\end{align*}
satisfying
\begin{align*}
W \, V^E_i = (V^C_i \oplus Y_i)\, W
\end{align*}
\begin{align*}
\tilde{W} \, V^{\tilde{E}}_i = (V^E_i \oplus \tilde{Y}_i) \,\tilde{W}
\end{align*}
We can define another unitary
\begin{align} \label{4.9}
Z := (W \otimes \eins_{\tilde{\CK}})\,\tilde{W}: \CH_{\tilde{E}}
\oplus (\Gamma \otimes \CD_{\tilde{E}}) \rightarrow \CH_C \oplus
(\Gamma \otimes \CD_C) \oplus \CK \oplus \tilde{\CK}
\end{align}
satisfying
\begin{align} \label{4.10}
Z \, V^{\tilde{E}}_i = (V^C_i \oplus Y_i \oplus \tilde{Y}_i) \,Z
\end{align}
Note further that $W, \tilde{W}$ and hence also $Z$ act identically on $\CH_C$.
By assumption the liftings $\uE$ of $\uC$ and $\tilde{\uE}$ of $\uE$ are reduced
and we have characteristic functions
\begin{align*}
M_{C,E} = P_{\Gamma \otimes \CD_C} W |_{\Gamma \otimes \CD_E}
\end{align*}
\begin{align*}
M_{E,\tilde{E}} = P_{\Gamma \otimes \CD_E} \tilde{W} |_{\Gamma \otimes \CD_{\tilde{E}}}
\end{align*}
They can be composed to yield a multi-analytic operator
\begin{align*}
M := M_{C,E} \, M_{E,\tilde{E}}: \;
\Gamma \otimes \CD_{\tilde{E}} \rightarrow \Gamma \otimes \CD_C
\end{align*}
Using \eqref{4.9} it is easily checked that
\begin{align*}
M = P_{\Gamma \otimes \CD_C} Z |_{\Gamma \otimes \CD_{\tilde{E}}}
\end{align*}
We conclude by \eqref{4.10} that the lifting $\tilde{\uE}$ of $\uC$
is associated to $M$ and hence, by Proposition 3.8, this lifting is reduced.
In fact, comparing with Definition 3.6, we see that $M$ is the characteristic
function, i.e., $M = M_{C,\tilde{E}}$.
\end{proof}

\section{Applications to Completely Positive Maps}

If $\uT = (T_1,\ldots,T_d)$ is a row contraction on a Hilbert space $\CK$ then we denote
by $\Phi_T$ the corresponding (normal) completely positive map on $\CB(\CK)$ given by
\begin{align} \label{5.1}
\Phi_T(\cdot) = \sum^d_{i=1} T_i \cdot T^*_i
\end{align}
If $d=\infty$ this should be understood as a SOT-limit. See for example \cite{Pau03} for the general
theory of completely positive maps, we shall only work with the concrete representation
\eqref{5.1}. The fact that $\uT$ is a row contraction implies that $\Phi_T(\eins) \leq \eins$,
i.e., $\Phi_T$ is contractive. It is unital ($\Phi_T(\eins) = \eins$) if and only if
$\uT$ is coisometric.

If $\uE$ is a contractive lifting of $\uC$ by $\uA$, i.e., $E_i =
\left(\begin{array}{cc}
C_i & 0   \\
B_i & A_i
\end{array}
\right)$
(for $i=1,\ldots,d$)
then an elementary computation shows that
\begin{align} \label{5.2}
\Phi_E
\left(\begin{array}{cc}
X_{11} & X_{12}   \\
X_{21} & X_{22}
\end{array}
\right)
\quad = \quad
\sum^d_{i=1}
\left(\begin{array}{cc}
\vspace{0.2cm}
C_i X_{11} C^*_i                    & C_i X_{11} B^*_i + C_i X_{12} A^*_i \\
                                    & B_i X_{11} B^*_i + B_i X_{12} A^*_i \\
B_i X_{11} C^*_i + A_i X_{21} C^*_i &+A_i X_{21} B^*_i + A_i X_{22} A^*_i
\end{array}
\right)
\end{align}
with $X_{11} \in \CB(\CH_C)$, $X_{12} \in \CB(\CH_A,\CH_C)$, $X_{21}
\in \CB(\CH_C,\CH_A)$, $X_{22} \in \CB(\CH_A)$. We denote by $p_C=
\left(\begin{array}{cc}
\eins & 0   \\
0 & 0
\end{array}
\right)$ and $p_A= \left(\begin{array}{cc}
0 & 0   \\
0 & \eins
\end{array}
\right)$
the orthogonal projections onto $\CH_C$ and $\CH_A$.
The following facts are immediate from \eqref{5.2}.

\begin{align} \label{5.3}
p_C \,(\Phi_E)^n
\left(\begin{array}{cc}
X & 0   \\
0 & 0
\end{array}
\right)
|_{\CH_C}
\quad = \quad
(\Phi_C)^n (X)
\end{align}
(for $n\in\N_0$ and $X\in\CB(\CH_C)$)

\begin{align} \label{5.4}
\Phi_E \left(\begin{array}{cc}
0 & 0   \\
0 & Y
\end{array}
\right)
\quad = \quad
\left(\begin{array}{cc}
0 & 0   \\
0 & \Phi_A(Y)
\end{array}
\right)
\end{align}
(for $Y\in \CB(\CH_A)$).
So $\Phi_E$ is a kind of (power) dilation of $\Phi_C$ \eqref{5.3} and an extension
of $\Phi_A$ \eqref{5.4}.

\begin{defn}
If $\CH_E = \CH_C \oplus \CH_A$, $\Phi_E: \CB(\CH_E) \rightarrow \CB(\CH_E)$,
$\Phi_C: \CB(\CH_C) \rightarrow \CB(\CH_C)$, $\Phi_A: \CB(\CH_A) \rightarrow \CB(\CH_A)$
are contractive normal completely positive maps such that \eqref{5.3} and \eqref{5.4} are
valid then we say that $\Phi_E$ is a contractive lifting of $\Phi_C$ by $\Phi_A$.
\end{defn}

We have seen that a contractive lifting of row contractions gives rise to a contractive
lifting of completely positive maps. The converse is also true: Let us assume \eqref{5.4}.
If $\Phi_E (\cdot) = \sum^d_{i=1} E_i \cdot E^*_i$ and we write
$E_i=
\left(\begin{array}{cc}
C_i & D_i   \\
B_i & A_i
\end{array}
\right)$
for the moment, then
\begin{align*}
\Phi_E
\left(\begin{array}{cc}
0 & 0   \\
0 & \eins
\end{array}
\right)
\quad = \quad
\left(\begin{array}{cc}
\vspace{0.2cm}
\sum^d_{i=1} D_i D^*_i & * \\
           *           & * \\
\end{array}
\right)
\end{align*}
and \eqref{5.4} implies that all the $D_i$ are zero, i.e., we have a lifting of row contractions.
So actually \eqref{5.4} implies \eqref{5.3} with some $\Phi_C$.

Note that if $\uE = \uV^C$, the mid of $\uC$, then $\Phi_E$ is a $*$-homomorphism and \eqref{5.3}
shows that the powers of $\Phi_E$ are a homomorphic dilation of the completely positive semigroup
formed by powers of $\Phi_C$. See \cite{BP94,Ar04,Go04} for further information about this kind of dilation theory.

The discussion above shows that we can use our theory of liftings for row contractions to study liftings of completely positive maps. If $\uE$ is a reduced lifting of $\uC$ by $\uA$ then we have
a characteristic function $M_{C,E}$. It is well known (see for example \cite{Pau03,Go04})
that in the decomposition
$\Phi_E {\cdot} = \sum^d_{i=1} E_i \cdot E^*_i$ the tuple $(E_1,\ldots,E_d)$ is not uniquely determined and that $\sum^d_{i=1} E'_i \cdot (E'_i)^*$ describes the same map if and only if
$\uE'$ is obtained from $\uE$ by multiplication with a unitary $d \times d$-matrix (with complex
entries). This does not change the characteristic function because the latter is defined as an intertwiner between objects which are transformed in the same way. Hence it is possible to think
of $M_{C,E}$ also as the characteristic function of a reduced lifting $\Phi_E$ of $\Phi_C$ by $\Phi_A$. (Of course we call this lifting reduced if the corresponding lifting of row contractions
is reduced.) Theorem 3.7 translates immediately into

\begin{cor}
Given $\Phi_C$, two reduced liftings $\Phi_E$ resp. $\Phi_{E'}$ of $\Phi_C$
by $\Phi_A$ resp. $\Phi_{A'}$ are conjugate, i.e.
\begin{align*}
\Phi_E = U^* \Phi_{E'} (UXU^*) U
\end{align*}
with a unitary $U: \CH_E \rightarrow \CH_{E'}$ such that $U |_{\CH_C} = \eins |_{\CH_C}$,
if and only if the corresponding characteristic functions are equivalent.
\end{cor}

\noindent
Corollary 5.2 generalizes Corollary 6.3 in \cite{DG07} where $dim\, \CH_C = 1$.
\\

In the following we confine ourselves mainly to liftings which are coisometric and subisometric
and give some concrete and useful results about the corresponding completely positive maps.

\begin{lem}
Let $\uE$ be a contractive lifting of a row contraction $\uC$ by a $*$-stable row contraction $\uA.$
Then for all $X_{12}, X_{21}, X_{22}$
\[
\Phi^n_E
\left(
\begin{array}{cc}
0  &  X_{12} \\
X_{21}  &  X_{22} \\
\end{array}
\right)
\to 0
\]
as $n \to \infty$ (SOT).
\end{lem}

\begin{proof}
$\Phi^n_E (p_A)$ decreases to zero in the strong operator topology
because of \eqref{5.4} and the assumption that $\uA$ is $*$-stable. Then also
$
\Phi^n_E
\left(
\begin{array}{cc}
0  &  0 \\
0  &  X_{22} \\
\end{array}
\right)
\to 0,
$
first for $0 \leq X_{22} \leq \|X_{22}\| \,p_A$, then for general $X_{22}$ by writing it as a linear combination
of positive elements. Using the Kadison-Schwarz inequality for completely positive maps
(cf. \cite{Ch74} or \cite{Pau03}, Chapter 3) we obtain
\begin{eqnarray*}
\Phi^n_E
\left(
\begin{array}{cc}
0  &  0 \\
X^*_{12}  &  0 \\
\end{array}
\right)\;
\Phi^n_E
\left(
\begin{array}{cc}
0  &  X_{12} \\
0  &  0 \\
\end{array}
\right)
&\leq & \Phi^n_E
\left(
\begin{array}{cc}
0  &  0 \\
0 & X^*_{12} X_{12} \\
\end{array}
\right)
\to 0
\end{eqnarray*}
and hence
$
\Phi^n_E
\left(
\begin{array}{cc}
0  &  X_{12} \\
0  &  0 \\
\end{array}
\right)
\to 0$.
Similarly
$
\Phi^n_E
\left(
\begin{array}{cc}
0  &  0 \\
X_{21} &  0 \\
\end{array}
\right)
\to 0$.
\end{proof}

\begin{thm}
Suppose the row coisometry $\uE$ is a lifting of $\uC$ by $\uA$.
Then the following assertions are equivalent:
\begin{itemize}
\item[(a)]
The lifting is subisometric.
\item[(b)]
$\uA$ is $*$-stable.
\item[(c)]
$(\Phi_E)^n (p_C) \to \eins \quad (n \to \infty, SOT)$
\item[(d)]
There is an order isomorphism between the fixed point sets of $\Phi_E$ and of $\Phi_C$
given by
\begin{align} \label{5.6}
\kappa: X \mapsto p_C \, X \, p_C
\end{align}
\end{itemize}
In this case, $\kappa$ is isometric on selfadjoint elements. If $x$ is a fixed point of $\Phi_C$
then we can reconstruct the preimage $\kappa^{-1}(x)$ as the SOT-limit
\begin{align} \label{5.7}
\lim_{n \to \infty}
(\Phi_E)^n
\left(
\begin{array}{cc}
x  &  0 \\
0  &  0 \\
\end{array}
\right)
\end{align}
\end{thm}

Recall further that by the results of Section 2 the liftings in Theorem 5.4 are
parametrized by $*$-stable row contractions $\uA$ with $dim\,\CD_{*,A} \leq dim\,\CD_C$
together with isometries $\gamma: \CD_{*,A} \rightarrow \CD_C$ and that they can be
explicitly constructed from these data. Theorem 5.4(d) tells us that (exactly) for such
liftings the maps $\Phi_E$ and $\Phi_C$ have closely related properties in terms of
their fixed points. We can identify this useful situation by checking the convenient conditions
(b) or (c).

\begin{proof}
By Proposition 2.3 a coisometric lifting $\uE$ of $\uC$ by $\uA$ is subisometric if and only if
$\uA$ is $*$-stable. Using \eqref{5.4} the latter means that
\begin{align*}
(\Phi_E)^n (p_A) \to 0 \quad (n \to \infty,\; SOT),
\end{align*}
which is equivalent to (c) because $\Phi_E$ is unital. Hence
$(a) \Leftrightarrow (b) \Leftrightarrow (c)$.

If $X = \left(\begin{array}{cc}
x & *   \\
* & *
\end{array}
\right)$
is a fixed point of $\Phi_E$ then it is immediate from \eqref{5.2} that $x$ is a
fixed point of $\Phi_C$. Hence $\kappa: X \mapsto p_C \, X \, p_C$ indeed maps fixed points
of $\Phi_E$ to fixed points of $\Phi_C$. (This is true for all contractive liftings.)
Now assume (a), i.e., the lifting is subisometric. Then
\begin{align*}
X = \Phi_E (X) = \lim_{n\to\infty} (\Phi_E)^n (X) = \lim_{n\to\infty} (\Phi_E)^n \left(
\begin{array}{cc}
x  &  0 \\
0  &  0 \\
\end{array}
\right),
\end{align*}
where the last equality follows from Lemma 5.3. Hence $\kappa$ is injective.

Let $\uV = (V_1,\ldots,V_d)$ simultaneously serve as
mid for $\uC$ and $\uE$. Then Theorem 5.1 in \cite{BJKW00} or Lemma 6.4 in the Appendix
of this paper show that
for every fixed point $x$ of $\Phi_C$ there exists $A'$ in the commutant
of $V_1,\ldots,V_d$ such that
$p_C A' p_C = x$. Define $X := p_E A' p_E$, where $p_E$ is the orthogonal projection
onto $\CH_E$.
Then, using the lifting property $E_i\,p_E = p_E \, V_i$ for $i=1,\ldots,d$ for the mid and
the fact that $\sum^d_{i=1} V_i V^*_i = \eins$ (because $\uE$ is coisometric also $\uV$ is
coisometric), we find that
\begin{align*}
\Phi_E (X) = \sum^d_{i=1} E_i X E^*_i = \sum^d_{i=1} E_i p_E A' p_E E^*_i
= \sum^d_{i=1} p_E V_i A' V^*_i p_E \\
= p_E A' \sum^d_{i=1} V_i V^*_i p_E
= p_E A' p_E = X
\end{align*}
So $X$ is a fixed point of $\Phi_E$ and clearly $\kappa(X)=x$. We conclude that $\kappa$ is
also surjective. The fact that $\kappa$ is isometric on selfadjoint elements is also a consequence of Lemma 6.4.

On the other hand, if the lifting $\uE$ of $\uC$ is not subisometric then the mid $\uV^C$
of $\uC$ is embedded on a proper reducing subspace $\hat{\CH}_C$ into the space $\hat{\CH}_E$
of the mid $\uV^E$ of $\uE$. Then $\eins_{\hat{\CH}_E}$ and $p_{\hat{\CH}_C}$ are two
different fixed points of $\Phi_{V^E}$. By Lemma 6.4 the map $\hat{X} \mapsto p_E \hat{X} p_E$
maps them into different fixed points of $\Phi_E$:
$p_E \eins_{\hat{\CH}_E} p_E = p_E \not= p_E p_{\hat{\CH}_C} p_E$.
If $\kappa: X \mapsto p_C \, X \, p_C$ from the fixed point set of $\Phi_E$ into
the fixed point set of $\Phi_C$ were injective then also
$p_C p_E p_C \not= p_C p_E p_{\hat{\CH}_C} p_E p_C$. But both sides are equal to $p_C$.
Hence in this case $\kappa$ is not injective. We have proved $(a) \Leftrightarrow (d)$.
\end{proof}

Recall that a unital completely positive map $\Phi_E$ is called {\it ergodic} if there are no other fixed
points than the multiples of the identity. By abuse of language we also call $\uE$ ergodic in this
case (as in \cite{DG07}).

\begin{prop}
Let $\uE$ be a coisometric lifting of $\uC$ by $\uA$. Then $\uE$ is
ergodic if and only if $\uC$ is ergodic and $\uA$ is $*$-stable.
\end{prop}

\begin{proof}
If $\uA$ is $*$-stable then use the equivalence $(b) \Leftrightarrow (d)$ in Theorem 5.4 and
infer from $\uC$ ergodic that also $\uE$ is ergodic.
Further note that, because $\uE$ is coisometric, we always have
\begin{align*}
\Phi_E
\left(
\begin{array}{cc}
\eins  &  0 \\
0  &  0 \\
\end{array}
\right)
\quad = \quad
\left(
\begin{array}{cc}
\eins  &  0 \\
0  &  \uB \uB^* \\
\end{array}
\right)
\quad \geq \quad
\left(
\begin{array}{cc}
\eins  &  0 \\
0  &  0 \\
\end{array}
\right)
\end{align*}
We say that $p_C \left(
\begin{array}{cc}
\eins  &  0 \\
0  &  0 \\
\end{array}
\right)$
is an increasing projection for $\Phi_E$.
Hence $(\Phi_E)^n (p_C)$ increases to a SOT-limit which clearly is a fixed point of $\Phi_E$.
\\
Now let $\uE$ be ergodic. Then all fixed points are multiples of
$\left(
\begin{array}{cc}
\eins  &  0 \\
0  &  \eins \\
\end{array}
\right)$
and because the left upper corner of
$(\Phi_E)^n (p_C)$ is always $\eins$
we have $(\Phi_E)^n (p_C) \to
\left(
\begin{array}{cc}
\eins  &  0 \\
0  &  \eins \\
\end{array}
\right)$.
We have verified Theorem 5.4(c) and now Theorem 5.4(d) and (b)
show that $\uC$ is ergodic and that $\uA$ is $*$-stable.
\end{proof}

This generalizes Proposition 2.3 in \cite{DG07} where $\CH_C$ is one dimensional
and hence $\uC$ ergodic is automatically fulfilled.
\\

The following provides an interesting example for the liftings considered above.
Let $\Phi_E: \CB(\CH_E) \rightarrow \CB(\CH_E)$ be any (normal) unital completely positive map
and let $\psi$ be a normal invariant state, i.e., $\psi \circ \Phi_E = \psi$.
Define $\CH_C$ to be the support of $\psi$ (cf.\,\cite{Ta01}) and let $\CH_A$ be the orthogonal complement,
so $\CH_E = \CH_C \oplus \CH_A$. Then $\uE = (E_1,\ldots,E_d)$ is a coisometric lifting
of $\uC = (C_1,\ldots,C_d)$ if we define $C^*_i := E^*_i |_{\CH_C}$ for $i=1,\ldots,d$.
In fact $p_C E^*_i p_C = E^*_i p_C$ for all $i$ by Lemma 6.1 of \cite{BJKW00}. Note that the
compression $\Phi_C$ has a faithful normal invariant state, the restriction $\psi_C$
of $\psi$ to $\CB(\CH_C)$. Conversely we can start with $\Phi_C$ and a faithful invariant state
$\psi_C$ and construct liftings $\Phi_E$ . They have normal invariant states given by
$\psi(X) := \psi_C (p_C X p_C)$. From Proposition 2.3 and Theorem 5.4 we conclude

\begin{cor}
Let $\Phi_C: \CB(\CH_C) \rightarrow \CB(\CH_C)$ be a (normal) unital completely positive map with a
faithful normal invariant state $\psi_C$. Then we have a one-to-one correspondence between
\begin{itemize}
\item[(a)]
(normal) unital completely positive maps $\Phi_E: \CB(\CH_E) \rightarrow \CB(\CH_E)$
with normal invariant state $\psi$ such that the support of $\psi$ is $\CH_C$ and
$\psi |_{\CB(\CH_C)} = \psi_C$, compression of $\Phi_E$ is $\Phi_C$ and
$(\Phi_E)^n (p_C) \to \eins \quad (n \to \infty,\, SOT)$
\item[(b)]
$*$-stable $\uA$ with $dim\,\CD_{*,A} \leq dim\,\CD_C$ together with isometries
$\gamma: \CD_{*,A} \rightarrow \CD_C$
\end{itemize}
There exist order isomorphisms $\kappa_E: X \mapsto p_C X p_C$ between the fixed point sets
of these maps $\Phi_E$ and the fixed point set of $\Phi_C$.
\end{cor}

In the special case when $\psi$ is an invariant vector state $\langle \xi, \cdot \,\xi \rangle$
of $\Phi_E$ we have the result that $\Phi_E$ is ergodic if and only if $(\Phi_E)^n (p_\xi) \to \eins
\; (n \to \infty,\, SOT)$, where $p_\xi$ is the orthogonal projection onto $\C \xi$, cf. \cite{Go04},
A.5.2.
Hence we obtain a classification of such maps. Here $\CD_C$ is
$(d-1)$-dimensional. This case has been further investigated in \cite{DG07}.

Corollary 5.6 is useful because many techniques only apply to completely positive maps with faithful invariant states, cf. \cite{Ku}.
It enables us to transfer information from the faithful to the non-faithful setting. For example, it is known that in the case of a faithful normal invariant state
the fixed point set is an algebra (cf. \cite{Ch74,FNW94,BJKW00}). Now $\kappa$ is an order isomorphism but it is not
in general multiplicative. In fact, there are examples of completely positive maps with a normal
invariant non-faithful state where the fixed point set is not an algebra
(cf. \cite{Ar69,Ar72,BJKW00}. If
Corollary 5.6 applies we can think of it as an (order isomorphic) deformation of an algebra.

\section{Appendix}

In Section 5 we needed a commutant lifting theorem (Theorem 5.1 of \cite{BJKW00}) which says that the fixed point set of a normal unital completely positive map is in one-to-one correspondence with the commutant of the Cuntz algebra representation generated by the mid. Below we give a variant
of the proof which is based on a Radon-Nikodym result for completely positive maps by W.Arveson.
This is a good way to think about it and it supports the understanding of the other arguments in the main text.

\begin{lem} \cite{Ar69}, Theorem 1.4.2 \\
If $\Psi$ is a completely positive map from a $C^*$-algebra $\CB$ to $\CB(\CH)$, with $\CH$ a Hilbert
space, then there exists an affine order isomorphism of the partially ordered set of operators
$\{ A' \in \pi(\CB)': 0 \leq A' \leq \eins \}$ onto $[0,\Psi]$. Here $\pi$ is the minimal Stinespring
representation of $\CB$ associated to $\Psi$ and $[0,\Psi]$ is the order interval containing all
completely positive maps $\Phi: \CB \rightarrow \CB(\CH)$ with $0 \leq \Phi \leq \Psi$. The order
relation for completely positive maps used here is $\Phi \leq \Psi$ if $\Psi - \Phi$ is completely
positive.
\\
Explicitly, if $\Psi(x) = W^* \pi(x) W$ is the minimal Stinespring representation of $\Psi$ then
$A' \in \pi(\CB)'$ corresponds to $\Phi = W^* A' \pi(x) W$.
\end{lem}

\begin{lem} \cite{BJKW00}, Corollary 2.4; \cite{Po03}, Theorem 2.1 \\
If $0 \leq D \leq \eins$ is a fixed point of the (normal unital completely positive) map
$\Phi_R(\cdot) = \sum^d_1 R_i \cdot R^*_i$ on $\CB(\CH)$ then there exists a completely
positive map $\Psi_D: \Od \rightarrow \CB(\CH), \; V_\alpha V^*_\beta \mapsto R_\alpha D R^*_\beta$. Here $\alpha,\beta \in \tilde{\Lambda}$ and $\Od$ is the Cuntz algebra generated by the $V_i$,
where $\uV = (V_1,\ldots,V_d)$ is a mid of $\uR = (R_1,\ldots,R_d)$.
\end{lem}

Using notation from the previous lemmas we get

\begin{lem}
There exists an affine order isomorphism $D \mapsto \Psi_D$ between
\[
\{ 0 \leq D \leq \eins: D \;\mbox{is a fixed point of}\; \Phi_R(\cdot) = \sum^d_1 R_i \cdot R^*_i
\;\mbox{on}\; \CB(\CH) \}
\]
and $[0,\Psi_\eins]$, where $\Psi_\eins$ is the completely positive map described in Lemma 6.2
with $D = \eins$, i.e., $\Psi_\eins: \Od \rightarrow \CB(\CH), \; V_\alpha V^*_\beta \mapsto R_\alpha R^*_\beta$.
\end{lem}

\begin{proof}
From $\Psi_\eins = \Psi_D + \Psi_{\eins - D}$ we see that $\Psi_D \in [0,\Psi_\eins]$
for all fixed points $ 0 \leq D \leq \eins$ of $\Phi_R$. On the other hand, if
$\Phi \in [0,\Psi_\eins]$ then by Lemma 6.1 with $\CB = \Od$ there exists $A' \in \pi(\CB)'$ with
$0 \leq A' \leq \eins$ such that $\Phi(x) = W^* A' \pi(x) W$, where $\Psi_\eins(x) = W^* \pi(x) W$
is a minimal Stinespring representation. Using that $(V_1,\ldots,V_d)$ is a mid of
$\uR = (R_1,\ldots,R_d)$ it is easily checked that $\Psi_\eins(x) = p \pi(x) p$ is such a minimal
Stimespring representation if $\pi$ is the Cuntz algebra representation generated by
$(V_1,\ldots,V_d)$ and $p$ is the projection onto the space $\CH$. (In $p \pi(x) p$ the $p$ on the
right hand side should be interpreted as the embedding of $\CH$ into the dilation space.)
\\
Hence if $x = V_\alpha V^*_\beta$ then we obtain
\[
\Phi(V_\alpha V^*_\beta) = p A' V_\alpha V^*_\beta p = p V_\alpha A' V^*_\beta p = p V_\alpha p A' p V^*_\beta p = R_\alpha p A' p R^*_\beta.
\]
We conclude that $\Phi = \Psi_D$ with $D := p A' p$. Clearly $0 \leq D \leq \eins$ and $D$ is a
fixed point of $\Phi_R$ (because $\uV$ is a coisometric lifting of $\uR$, i.e.,
$\sum^d_{i=1} V_i V^*_i = \eins$ and $R_i p = p V_i$ for all $i$).
The correspondence
is bijective ($\Psi_D(\eins) = D$) and it clearly respects the order.
\end{proof}

\begin{lem} \cite{BJKW00}, Theorem 5.1 \\
There is an affine order isomorphism between
$\{ 0 \leq D \leq \eins: D \;\mbox{is a fixed point}\\
\mbox{of}\;
 \Phi_R(\cdot) = \sum^d_1 R_i \cdot R^*_i
\;\mbox{on}\; \CB(\CH) \}$ and $\{ A' \in \pi(\Od)': 0 \leq A' \leq \eins \}$, where $\pi$ is the
Cuntz algebra representation generated by the mid
$\uV = (V_1,\ldots,V_d)$ of $\uR = (R_1,\ldots,R_d)$. It is given by $A' \mapsto p A' p$,
where $p$ is the projection onto the space $\CH$. The isomorphism is isometric on the
selfadjoint parts.
\end{lem}

\begin{proof}
For the first part we only have to add to the arguments in the proof of Lemma 6.3 the
reminder that by Lemma 6.1 the correspondence between $\{ A' \in \pi(\Od)': 0 \leq A' \leq \eins \}$
and $[0, \Psi_\eins]$ is a bijection. As pointed out in \cite{BJKW00}, Section 4, it is isometric on the
selfadjoint parts because $\eins$ is mapped to $\eins$ (identities on different Hilbert spaces) and for selfadjoint elements $y$ we have
$\|y\| = \inf \{ \alpha > 0: -\alpha \eins \leq y \leq \alpha \eins \}$.
\end{proof}

\vspace{1cm}
\noindent Institut f\"ur Mathematik und Informatik,
E-M-A-Universit\"at Greifswald, F-L-Jahnstr. 15a, D-17487
Greifswald, Germany

\vspace{.5cm}

\noindent Department of Mathematics, University of Reading,
Whiteknights, P.O.Box 220, Reading, RG6 6AX, UK

\end{document}